%% file: agt-5-39.tex
\documentclass{gtart_h}
\input agtout

\lognumber{39}
\volumenumber{5}
\volumeyear{2005}
\papernumber{39}
\pagenumbers{923}{964}
\received{16 February 2005} 
\accepted{7 July 2005}
\published{5 August 2005}

\usepackage{amssymb}
\usepackage[matrix,arrow,curve]{xy}
\usepackage{epsfig}
\usepackage{psfrag}
\usepackage{amsmath}

\def\figref#1{\hyperlink{#1anchor}{Figure~\ref*{#1}}}
\def\anchor#1{\noindent\hypertarget{#1anchor}{\smash{$\phantom{99}$}}}

\numberwithin{equation}{section}
\numberwithin{figure}{section}

\newcommand{\ccoeg}{:=}

\newtheorem{thm}{Theorem}[section]

\newtheorem{Theorem}[thm]{Theorem}

\newtheorem{Lemma}[thm]{Lemma}

\newtheorem{Proposition}[thm]{Proposition}
\newtheorem{Corollary}[thm]{Corollary}

\theoremstyle{remark}
\newtheorem{Definition}[thm]{Definition}

\newtheorem{Remark}[thm]{Remark}

\newtheorem{Example}[thm]{Example}
\newtheorem{Examples}[thm]{Examples}

\theoremstyle{plain}
\swapnumbers
\newtheorem{cote}[subsubsection]{}
\newenvironment{ccote}{\cote}{}

\newtheorem{qote}[thm]{}
\newenvironment{cqote}{\qote}{}

\renewcommand{\mod}{{\rm mod\,}}
\newcommand{\ri}{\kappa}

\newcommand{\bbd}{{\mathbb{D}}}
\newcommand{\bbs}{{\mathbb{S}}}
\newcommand{\bbp}{{\mathbb{P}}}
\newcommand{\bbe}{{\mathbb{E}}}
\newcommand{\bbf}{{\mathbb{F}}}
\newcommand{\bbr}{{\mathbb{R}}}
\newcommand{\bbc}{{\mathbb{C}}}

\newcommand{\bbz}{{\mathbb{Z}}}
\newcommand{\bbq}{{\mathbb{Q}}}
\newcommand{\bbn}{{\mathbb{N}}}

\newcommand{\calb}{{\mathcal B}}
\newcommand{\calc}{{\mathcal C}}

\newcommand{\cali}{{\mathcal I}}

\newcommand{\calk}{{\mathcal K}}

\newcommand{\calo}{{\mathcal O}}

\newcommand{\calt}{{\mathcal T}}
\newcommand{\calu}{{\mathcal U}}
\newcommand{\calv}{{\mathcal V}}

\newcommand{\pcirc}{\kern .7pt {\scriptstyle \circ} \kern 1pt}
\newcommand{\iso}{\approx}
\newcommand{\hfl}[1]{\buildrel{#1}\over{\longrightarrow}}
\newcounter{exo}

\newcommand{\mun}{{-1}}

\newcommand{\llangle}[2]{\langle #1 ,#2 \rangle}

\newcommand{\onto}{\to\kern-6.5pt\to}
\newcommand{\into}{\hookrightarrow}
\newcommand{\algt}{\mathfrak{t}}
\newcommand{\algg}{\mathfrak{g}}
\newcommand{\algh}{\mathfrak{h}}
\newcommand{\algl}{\mathfrak{l}}
\newcommand{\proref}[1]{Proposition~\ref{#1}}
\newcommand{\remref}[1]{Remark~\ref{#1}}
\newcommand{\lemref}[1]{Lemma~\ref{#1}}

\newcommand{\thref}[1]{Theorem~\ref{#1}}

\newcommand{\dunion}{\,\hbox{\raisebox{.48ex}{$\scriptstyle \coprod\,$}}}
\newcommand{\rev}{R^{ev}}

\newcommand{\lt}[1]{\ell t_{#1}}
\newcommand{\csp}{conjugation space}
\newcommand{\cpa}{conjugation pair}
\newcommand{\hfra}{$H^*$-frame}
\newcommand{\scr}{\scriptscriptstyle}

\newcommand{\ef}{equivariantly formal}

\newcommand{\scc}{spherical conjugation complex}
\newcommand{\taub}{$\tau$-bundle}
\newcommand{\tho}{{\rm Thom}}
\newcommand{\ceb}{conjugate equivariant bundle}

\title{Conjugation spaces}
\author{Jean-Claude Hausmann\\Tara Holm\\Volker Puppe}
\asciiauthors{Jean-Claude Hausmann, Tara Holm and Volker Puppe}
\shortauthors{Hausmann, Holm and Puppe}
\address{Section de math\'ematiques, 2-4, rue du Li\`evre\\CP 64
CH-1211 Gen\`eve 4, Switzerland\\\smallskip\\Department of Mathematics,
University of Connecticut\\Storrs CT 06269-3009, 
USA\\\smallskip\\Universit\"at Konstanz, 
Fakult\"at f\"ur Mathematik\\Fach D202, D-78457 Konstanz, Germany}

\asciiaddress{Section de mathematiques, 2-4, rue du Lievre\\CP 64
CH-1211 Geneve 4, Switzerland\\Department of Mathematics, University
of Connecticut\\Storrs CT 06269-3009, USA\\Universitat Konstanz,
Fakultat fur Mathematik\\Fach D202, D-78457 Konstanz, Germany}

\asciiemail{hausmann@math.unige.ch, tsh@math.uconn.edu, Volker.Puppe@uni-konstanz.de}
\gtemail{\mailto{hausmann@math.unige.ch}, \mailto{tsh@math.uconn.edu}, 
\mailto{Volker.Puppe@uni-konstanz.de}}

\begin{abstract}
There are classical examples of spaces $X$ with an involution $\tau$
whose ${\rm mod} 2$-comhomology ring resembles that of their
fixed point set $X^\tau$: there is a ring isomorphism
$\kappa\co H^{2*}(X)\approx H^{*}(X^\tau)$.
Such examples include complex Grassmannians, toric manifolds,
polygon spaces. In this paper,
we show that the ring isomorphism $\kappa$ is part of an
interesting structure in equivariant cohomology called an {\it $H^*$-frame}.
An $H^*$-frame, if it exists, is natural and unique. A space with
involution admitting an $H^*$-frame is called a {\it\csp}.
Many examples of conjugation spaces are constructed, for instance
by successive adjunctions of cells homeomorphic to a disk in $\bbc^k$
with the complex conjugation. A
compact symplectic manifold, with an anti-symplectic involution
compatible with a Hamiltonian action of a torus $T$,
is a \csp, provided $X^T$ is itself a \csp.
This includes the co-adjoint orbits of any semi-simple compact Lie group,
equipped with the Chevalley involution.
We also study conjugate-equivariant complex vector bundles (``real bundles''
in the sense of Atiyah) over a \csp\ and show that the isomorphism $\kappa$
maps the Chern classes onto the Stiefel-Whitney classes of
the fixed bundle.
\end{abstract}

\asciiabstract{There are classical examples of spaces X with an
involution tau whose mod 2-comhomology ring resembles that of their
fixed point set X^tau: there is a ring isomorphism kappa: H^2*(X) -->
H^*(X^tau).  Such examples include complex Grassmannians, toric
manifolds, polygon spaces. In this paper, we show that the ring
isomorphism kappa is part of an interesting structure in equivariant
cohomology called an H^*-frame.  An H^*-frame, if it exists, is
natural and unique. A space with involution admitting an H^*-frame is
called a conjugation space.  Many examples of conjugation spaces are
constructed, for instance by successive adjunctions of cells
homeomorphic to a disk in C^k with the complex conjugation. A compact
symplectic manifold, with an anti-symplectic involution compatible
with a Hamiltonian action of a torus T, is a conjugation space,
provided X^T is itself a conjugation space.  This includes the
co-adjoint orbits of any semi-simple compact Lie group, equipped with
the Chevalley involution.  We also study conjugate-equivariant complex
vector bundles (`real bundles' in the sense of Atiyah) over a
conjugation space and show that the isomorphism kappa maps the Chern
classes onto the Stiefel-Whitney classes of the fixed bundle.}


\begin{document}
\maketitle

%\tableofcontents

\section{Introduction}\label{intro}

In this article, we study topological spaces equipped with a
continuous involution. We are motivated by the example of the
complex Grassmannian ${\rm Gr}(k,\bbc^n)$ of complex $k$-vector
subspaces of $\bbc^n$ ($n\leq\infty$), with the involution complex
conjugation. The fixed point set of this involution is the real
Grassmannian ${\rm Gr}(k,\bbr^n)$. It is well known that there is
a ring isomorphism $\kappa\co  H^{2*}({\rm Gr}(k,\bbc^n)) \approx
H^{*}({\rm Gr}(k,\bbr^n))$ in cohomology (with
$\bbz_2$-coefficients) dividing the degree of a class in half.
Other such isomorphisms have been found for natural
involutions on smooth toric manifolds \cite{DJ} and polygon spaces
\cite[\S\,9]{HK}. The
significance of this isomorphism property was first discussed
by A.~Borel and A.~Haefliger \cite{BH} in
the framework of analytic geometry. For more recent consideration
in the context of real algebraic varieties, see \cite{vH1,vH2}.

The goal of this paper is to show that, for the above examples and
many others, the ring isomorphism $\kappa$ is part of an
interesting structure in equivariant cohomology.  We will use
$H^*$ to denote singular cohomology, taken with $\bbz_2$
coefficients.  For a group $C$, we will let $H_C^*$ denote
$C$-equivariant cohomology with $\bbz_2$ coefficients, using the
Borel construction \cite{Bo}. Let $\tau$ be a (continuous) involution on
a topological space $X$. We view this as an action of the cyclic
group of order two, $C=\{I,\tau\}$. Let $\rho\co H^{2*}_C(X)\to
H^{2*}(X)$ and $r\co H^{*}_C(X)\to H^{*}_C(X^\tau)$ be the
restriction homomorphisms in cohomology. We use that
$H^*_C(X^\tau)=H^*(X^\tau\times BC)$ is naturally isomorphic to
the polynomial ring $H^*(X^\tau)[u]$ where $u$ is of degree one.
Suppose that $H^{odd}(X)=0$. A {\it cohomology frame} or {\it
\hfra}\  for $(X,\tau)$ is a pair $(\kappa,\sigma)$, where
\renewcommand{\labelenumi}{(\alph{enumi})}
\begin{enumerate}
\item
$\kappa\co  H^{2*}(X)\to H^{*}(X^\tau)$ is an additive isomorphism
dividing the degrees in half; and
\item $\sigma\co  H^{2*}(X)\to H^{2*}_C(X)$ is
an additive section of $\rho$.
\end{enumerate}
In addition, $\kappa$ and $\sigma$ must satisfy the {\it
conjugation equation}
\begin{equation}\label{defalceq-intro}
r\pcirc\sigma(a) = \kappa(a)u^m + \lt{m},
\end{equation}
for all $a\in H^{2m}(X)$ and all $m\in\bbn$, where $\lt{m}$
denotes some polynomial in the variable $u$ of degree less than $m$.

An involution admitting a \hfra\ is called a {\it conjugation} and
a space together with a conjugation is called a {\it \csp}.
Required  to be only additive maps, $\kappa$ and $\sigma$ are
often easy to construct degree by degree. But we will show in the
``multiplicativity theorem'' in Section~\ref{alc} that in fact
$\sigma$ and $\kappa$ are ring homomorphisms.  Moreover,
given a $C$-equivariant map $f\co Y\to X$ between spaces with
involution, along with \hfra s $(\sigma_X,\kappa_X)$ and
$(\sigma_Y,\kappa_Y)$, we have $H^*_Cf \pcirc\sigma_X
=\sigma_Y\pcirc H^*f$ \ and \ $H^*f^\tau \pcirc\kappa_X
=\kappa_Y\pcirc H^*f$. In particular, the \hfra\ for a conjugation
is unique.

As an example of a \csp, one has the complex projective space $\bbc P^k$
($k\leq\infty$),
with the complex conjugation as involution.
If $a$ is the generator of
$H^2(\bbc P^k)$ and $b=\kappa(a)$ that of $H^1(\bbr P^k)$,
we will see that
the conjugation equation has the form
$r\pcirc\sigma(a^m)=(bu+b^2)^m$ (Example~\ref{excpn}).

The complex projective spaces are particular cases of 
\scc es, which constitute our main class of examples.  
A {\it \scc} is a space (with involution)
obtained from the empty set by countably many successive adjunction 
of collections of conjugation cells.
A {\it conjugation cell} (of dimension $2k$)
is a space with involution which is equivariantly homeomorphic to 
the closed disk of radius $1$ in $\bbr^{2k}$, equipped with a linear involution 
with exactly $k$ eigenvalues equal to $-1$. At each step, the collection of
conjugation cells consists of cells of the same dimension but,
as in \cite{Gr}, the adjective ``spherical'' is a warning that
these dimensions do not need to be increasing.
We prove that every \scc\ is a \csp. 
There are many examples of these; for instance, there are infinitely many
$C$-equivariant homotopy types of \scc es with three conjugation
cells, one each in dimensions $0$, $2$ and $4$. We prove that for
a $C$-equivariant fibration (with a compact Lie group as structure group)
whose fiber is a \csp\ and whose base is a \scc, then its total space
is a \csp. 

Schubert cells for Grassmannians are conjugation cells, so these spaces
are \scc es and therefore \csp s. This generalizes
in the following way.  Let $X$ be
a space together with an involution $\tau$ and a continuous action
of a torus $T$. We say that $\tau$ is {\it compatible} with this
torus action if $\tau(g\cdot x)=g^\mun\cdot\tau(x)$ for all $g\in
T$ and $x\in X$.  It follows that $\tau$ induces an involution on
the fixed point set $X^T$ and an action of the 
$2$-torus $T_2$ (the elements of order $2$) of $T$ on $X^\tau$.
We are particularly interested in
the case when $X$ is a compact symplectic manifold for which the
torus action is Hamiltonian and the compatible involution is
smooth and anti-symplectic. Using a Morse-Bott function obtained
from the moment map for the $T$-action, we prove that if $X^T$ is
a \csp\ (respectively a \scc), then $X$ is a \csp\
(respectively a \scc).  In addition, we prove that
the involution induced on the Borel construction $X_T$ is a conjugation.
The relevant isomorphism $\bar\kappa$ takes the form of a natural
ring isomorphism
$$
\bar\kappa\co  H^{2*}_T(X)\hfl{\approx}{} H^{*}_{T_2}(X^\tau).
$$
Examples of such Hamiltonian spaces
include co-adjoint orbits of any semi-simple compact Lie group, with
the Chevalley involution, smooth toric manifolds and polygon spaces.
Consequently, these examples are \scc es.
For the co-adjoint orbits of $SU(n)$ this was proved 
earlier by C.~Schmid \cite{Sc} and D. Biss, V.~Guillemin and the
second author \cite{BGH}.
The category of \csp s is closed under various
operations, including direct products, connected sums and,
under some hypothesis, under symplectic reduction
(generalizing \cite{GH}; see Subsection~\ref{sred}).
This yields more examples of \csp s.

Over spaces with involution, it is natural to study 
{\it \ceb s}, identical to the
``real bundles'' introduced by Atiyah \cite{At}. These are complex
vector bundles $\eta=(E\hfl{p}{} X)$ together with an involution
$\hat\tau$ on $E$, which covers $\tau$ and is conjugate linear on
each fiber. Then $E^{\hat\tau}$ is a real bundle $\eta^\tau$ over
$X^\tau$. In Section~\ref{tabdles}, we prove several results on
\ceb s, among them that if $\eta=(E\hfl{p}{} X)$ is a \ceb\ over a
\csp, then the Thom space is a \csp. These results are used in the proof of
the aforementioned theorems in symplectic geometry. 
Finally, when the basis of a \ceb\ is a \scc, we prove that
$\kappa(c(\eta))=w(\eta^{\tau})$, where $c()$ denotes the ($\mod\
2$) total Chern class and $w()$ the total Stiefel-Whitney class.

\medskip
{\bf Acknowledgments}\qua
Anton Alekseev gave us precious
suggestions for Subsection~\ref{Cheva} on the Chevalley involution
on coadjoint orbits. Conversations with Matthias Franz were very helpful.
The first two authors are also grateful to Sue Tolman for pointing out
a gap in an earlier stage of this project. Finally, we thank 
Martin Olbermann for useful observations.

The three authors thank
the Swiss National Funds for Scientific Research for its support.
The second author was supported in part by a National Science Foundation 
Postdoctoral Fellowship.

%%%%%%%%%%%%%%%%%%%%%%%%%%%%%%%%%%%%%%%%%%%%%%%%%%%%%%%%%%%
\section{Preliminaries}\label{Sdefi}

Let $\tau$ be a (continuous) involution on a space $X$.
This gives rise to a continuous action of the cyclic
group $C=\{1,\tau\}$ of order $2$. 
The {\it real locus} $X^\tau\subset X$ is the subspace
of $X$ formed by the elements which are fixed by $\tau$.

Unless otherwise specified, all the cohomology groups are taken with
$\bbz_2$-coeffi\-cients. A pair $(X,Y)$ is an {\it even cohomology
pair} if $H^{odd}(X,Y)=0$; a space $X$ is an {\it even cohomology
space} if $(X,\emptyset)$ is an even cohomology pair.

\begin{cqote}\label{polu} \rm
Let $R$ be the graded ring $R=H^*(BC)=H^*_C(pt)=\bbz_2[u]$,
where $u$ is in degree $1$. We denote by $\rev$ the subring of
$R$ of elements of even degree.

As $C$ acts trivially on the real locus $X^\tau$,
there is a natural identification
$EC\times_C X^\tau\stackrel{\approx}\to BC\times X^\tau$. The K\"unneth formula
provides a ring isomorphism\\
 $R\otimes H^*(X^\tau,Y^\tau)\stackrel{\approx}\to
H^*_C(X^\tau,Y^\tau)$ and $R\otimes H^*(X^\tau,Y^\tau)$ is naturally isomorphic to the
polynomial ring $H^*(X^\tau,Y^\tau)[u]$. We shall thus often use the ``K\"unneth isomorphism''
$K\co H^*(X^\tau,Y^\tau)[u]\hfl{\approx}{} H^*_C(X^\tau,Y^\tau)$
to identify these two rings. The naturality of $K$ gives the following:

\begin{Lemma}\label{hfpri}
Let $f\co (X_2,Y_2)\to (X_1,Y_1)$ be a continuous $C$-equivariant map
between pairs with involution. Let
$f^\tau\co (X_2^\tau,Y_2^\tau)\to (X_1^\tau,Y_1^\tau)$ be the restriction
of $f$ to the fixed point sets. Then, the following diagram
$$\xymatrix@C+3pt@M+4pt@R-4pt{%
H^*(X_1^\tau,Y_1^\tau)[u] \ar[d]_(0.50){K}^(0.50){\approx}
\ar[r]^(0.50){H^*f^\tau[u]} &
H^*(X_2^\tau,Y_2^\tau)[u] \ar[d]^(0.50){\approx}_(0.50){K}  \\
%%%%%% ROW 2
H^*_C(X_1^\tau,Y_1^\tau) \ar[r]^(0.50){H^*_Cf^\tau}  &
H^*_C(X_2^\tau,Y_2^\tau)
}$$
is commutative, where $H^*f^\tau[u]$
is the polynomial extension of $H^*f^\tau$. \qed
\end{Lemma}
\end{cqote}

\begin{cqote}\label{eqform} Equivariant formality\qua \rm
Let $X$ be a space with an involution $\tau$ and let $Y$ be a
$\tau$-invariant subspace of $X$ (i.e.\ $\tau(Y)=Y$). Following
\cite{GKM}, we say that the pair $(X,Y)$ is {\it \ef\ } (over
$\bbz_2$) if the map $(X,Y)\to (EC\times_C X,EC\times_C Y)$ is
totally non-homologous to zero. That is, the restriction
homomorphism $\rho\co H^*_C(X,Y)\to H^*(X,Y)$ is surjective. A space
$X$ with involution is {\it \ef\ } if the pair $(X,\emptyset)$ is
\ef.

If $(X,Y)$ is \ef, one can choose, for each $k\in\bbn$,
a $\bbz_2$-linear map $\sigma\co  H^k(X,Y)\to H^k_C(X,Y)$
such that $\rho\pcirc\sigma={\rm id}$. This gives an additive section
$\sigma\co  H^*(X,Y)\to H^*_C(X,Y)$ of $\rho$ which gives rise
to a map
\begin{equation}\label{LH}
\hat\sigma\co  H^*(X,Y)[u] \to H^*_C(X,Y).
\end{equation}
As in \ref{polu}, we use the ring isomorphism $H^*(X,Y)\otimes R
\hfl{\approx}{} H^*(X,Y)[u]$. As $H^*(BC)=R$, the Leray-Hirsch
theorem (see e.g.\ \cite[Theorem 5.10]{McC}) then implies that
$\hat\sigma$ is an isomorphism of $R$-modules. But $\hat\sigma$ is
not in general  an isomorphism of rings. This is the case if and
only if the section $\sigma$ is a ring homomorphism but such 
ring-sections do not usually exist.
\end{cqote}

%%%%%%%%%%%%%%%%%%%%%%%%%%%%%%%%%%%%%%%%%%%%%%%%%%%%%%%%%%%%
\section{Conjugation pairs and spaces}\label{alc}

\subsection{Definitions and the multiplicativity theorem}\label{ssalc}
Let $\tau$ be an involution on a space $X$ and
let $Y$ be a $\tau$-invariant subspace of $X$.
Let $\rho\co  H^{2*}_C(X,Y)\to H^{2*}(X,Y)$ and
$r\co H^{*}_C(X,Y)\to H^{*}_C(X^\tau,Y^\tau)$ be the restriction homomorphisms.
A {\it cohomology frame} or {\it \hfra}\
for $(X,Y)$ is a pair $(\kappa,\sigma)$, where
\renewcommand{\labelenumi}{(\alph{enumi})}
\begin{enumerate}
\item $\kappa\co  H^{2*}(X,Y)\to H^{*}(X^\tau,Y^\tau)$ is
an additive isomorphism dividing the degrees in half; and
\item $\sigma\co  H^{2*}(X,Y)\to H^{2*}_C(X,Y)$ is
an additive section of $\rho$.
\end{enumerate}
Moreover, $\kappa$ and $\sigma$ must satisfy the {\it conjugation equation}
\begin{equation}\label{defalceq}
r\pcirc\sigma(a) = \kappa(a)u^m + \lt{m}
\end{equation}
for all $a\in H^{2m}(X)$ and all $m\in\bbn$, where $\lt{m}$
denotes any polynomial in the variable $u$ of degree less than $m$.

An involution admitting a \hfra\ is called a {\it conjugation}.
An even cohomology pair together with a conjugation is called a {\it \cpa}.
An even-cohomology space $X$ together with an involution
is a {\it\csp\ } if the pair $(X,\emptyset)$ is a \cpa .
Observe that the existence of $\sigma$ is equivalent to $(X,Y)$
being \ef. We shall see in Corollary \ref{alcstrunique} that the \hfra\ for
a conjugation is unique.

\begin{Remark}\label{kappaHo} The map $\kappa$ coincides on $H^0(X,Y)$ with the
restriction homomorphism $\tilde r\co H^0(X,Y)\to H^0(X^\tau,Y^\tau)$. \ \rm
Indeed, the following diagram
$$\xymatrix@C-3pt@M+2pt@R-4pt{%
H^0_C(X,Y) \ar[d]_(0.50){r} \ar[rr]^(0.50){\rho}_(0.50){\approx}
&& H^0(X,Y) \ar[d]^(0.50){\tilde r}
\\ %%%%%% ROW 2
H^0_C(X^\tau,Y^\tau) \ar[r]^(0.43){=}  &
H^0(X^\tau,Y^\tau)[u]^{[0]} \ar[r]^(0.55){=} &
H^0(X^\tau,Y^\tau)
}$$
is commutative. Therefore, using Equation~\eqref{defalceq},
one has for $a\in H^0(X,Y)$ that
$\kappa(a) = r\pcirc\sigma(a) = \tilde r (a)$. As a consequence,
if $X$ is a \csp, then $\pi_0(X^\tau)\approx \pi_0(X)$.
This implies that $\tau$ preserves each path-connected component of $X$.
\end{Remark}

\begin{Remark}\label{reduce}
Let $X$ be an path-connected space with an involution $\tau$.
Suppose that $X^\tau$ is non-empty and path-connected.
Let $pt\in X^\tau$. Then, $X$ is a \csp\ if and only if
$(X,pt)$ is a \cpa.
\end{Remark}

The remainder of this section is devoted to
establishing the fundamental facts about \cpa s and spaces,
and providing several important examples.

\begin{Theorem}[The multiplicativity theorem]\label{sigmaring}
Let $(\kappa,\sigma)$ be a \hfra\ for a conjugation $\tau$
on a pair $(X,Y)$.
Then $\kappa$ and $\sigma$ are ring homomorphisms.
\end{Theorem}

%\preu
\begin{proof}
We first prove that
\begin{equation}\label{sigmaring-eq1}
\sigma(ab)=\sigma(a)\sigma(b)
\end{equation}
for all $a\in H^{2k}(X,Y)$ and $b\in H^{2l}(X,Y)$. Let $m=k+l$.
Since $\rho\co  H^0_C(X,Y)\to H^0(X,Y)$ is an isomorphism,
Equation~\eqref{sigmaring-eq1} holds for $m=0$,
and thus we may assume that $m>0$.
As one has  $\rho\pcirc\sigma(ab)=\rho(\sigma(a)\sigma(b))$,
Equation~\eqref{sigmaring-eq1} holds true
modulo $\ker\rho$ which is the ideal generated by $u$.
As $H^*(X,Y)$ is concentrated
in even degrees, this means that
\begin{equation}\label{sigmaring-eq2}
\sigma(ab)=\sigma(a)\sigma(b)
+ \sigma(d_{2m-2})u^2 + \cdots + \sigma(d_{0})u^{2m} \ ,
\end{equation}
with $d_i\in H^{i}(X,Y)$. We must prove that $d_{2m-2}=\cdots=d_{0}=0$.

Let us apply $r\co H^{*}_C(X,Y)\to H^{*}_C(X^\tau,Y^\tau)$ to
Equation~\eqref{sigmaring-eq2}. The left hand side gives
\begin{equation}\label{sigmaring-eq4}
r\pcirc\sigma(ab)=
\kappa(ab)u^m + \lt{m}
\end{equation}
while the right hand side gives
\begin{equation}\label{sigmaring-eq3}
r\pcirc\sigma(ab)=
\kappa(a)\kappa(b)u^m + \lt{m} + (\kappa(d_{2m-2})u^{m-1} + \lt{m-1} )u^2
+ \cdots + \kappa(d_0)u^{2m} \, .
\end{equation}
Equations \eqref{sigmaring-eq4} and \eqref{sigmaring-eq3}
imply that
\begin{equation}\label{sigmaring-eq5}
r\pcirc\sigma(ab)= \kappa(d_0)u^{2m} + \lt{2m} \, .
\end{equation}
Comparing Equations~\eqref{sigmaring-eq4}
and~\eqref{sigmaring-eq5}, we deduce that $d_0=0$,
since $\kappa$ is injective. 
Then Equation~\eqref{sigmaring-eq2} implies that
\begin{equation}\label{sigmaring-eq5b}
r\pcirc\sigma(ab)= \kappa(d_2)u^{2m-1} + \lt{2m-1} \, .
\end{equation}
Again, comparing Equations~\eqref{sigmaring-eq4}
and~\eqref{sigmaring-eq5b}, we deduce that $d_2=0$.
This process continues until $d_{2m-2}$,
showing that each $d_i$ vanishes in
Equation~\eqref{sigmaring-eq2}, which proves that
$\sigma(ab)=\sigma(a)\sigma(b)$.

To establish that $\kappa(ab)=\kappa(a)\kappa(b)$ for
$a,b$ as above, we use the fact that
$r\pcirc\sigma(ab)=r\pcirc\sigma(a)\cdot r\pcirc\sigma(b)$
together with Equation~\eqref{defalceq} to conclude that
$$
\kappa(ab)u^m + \lt{m} = (\kappa(a)u^k + \lt{k})\,(\kappa(b)u^l + \lt{l})
= \kappa(a)\kappa(b)u^m + \lt{m} \ .
$$
Therefore, $\kappa$ is multiplicative. 
\end{proof}

By the Leray-Hirsch theorem, the section $\sigma$ gives rise
to an isomorphism of $R$-modules
$$\hat\sigma\co  H^{*}(X,Y)[u]\stackrel{\approx}\to H^{*}_C(X,Y)$$
(see \eqref{eqform}). As $\sigma$ is a ring homomorphism by
Theorem~\ref{sigmaring}, one has the following corollary,
which completely computes the ring $H^*_C(X,Y)$
in terms of $H^*(X,Y)$.

\begin{Corollary}\label{sigmaRalgebrabis}
Let $(\kappa,\sigma)$ be a \hfra\ for a conjugation
on a pair $(X;Y)$.
Then $\hat\sigma\co  H^{*}(X,Y)[u]\stackrel{\approx}\to H^{*}_C(X,Y)$
is an isomorphism of $R$-algebras. \qed 
\end{Corollary}

Finally, there is a unique map
$\hat\kappa\co  H^{2*}_C(X,Y)\to H^{*}_C(X^\tau,Y^\tau)$ such that the following diagram
\begin{equation}\label{defhatkappa}
\xymatrix@C-3pt@M+2pt@R-4pt{%
H^{2*}(X,Y)\otimes R^{ev} \ar[d]_{\kappa\otimes\alpha}
\ar[r]^(0.54){\hat\sigma}_(0.54){\approx} &  H^{2*}_C(X,Y)
\ar[d]_{\hat\kappa}
\\
H^*(X^\tau,Y^\tau)\otimes R \ar[r]^(0.54){K}_(0.54){\approx} &
H^{*}_C(X^\tau,Y^\tau) }
\end{equation}
is commutative, where $K$ comes from the K\"unneth formula.
The map $\hat\kappa$ is an isomorphism of
$(\rev,R)$-algebras over $\alpha\co \rev\to R$.

We now turn to examples of conjugation spaces and pairs.

\begin{Example}[Conjugation cells]\label{disk}\rm
Let $D=D^{2k}$ be the closed disk of radius $1$ in $\bbr^{2k}$, equipped
with an involution $\tau$ which is topologically conjugate to a
linear involution with exactly~$k$ eigenvalues equal to~$-1$.
We call such a disk a {\it conjugation cell} of dimension~$2k$.
Let $S$ be the boundary of $D$.
The fixed point set is then homeomorphic to a disk of dimension~$k$.

As $H^*(D,S)$ is concentrated in degree $2k$, the restriction
homomorphism\\ \vskip -3.3mm $\rho\co H^{2k}_C(D,S)\to H^{2k}(D,S)$
is an isomorphism.
Set $\sigma=\rho^\mun\co H^{2k}(D,S)\to H^{2k}_C(D,S)$.
This shows that $(D,S)$ is \ef.
The cohomology $H^*(D^\tau,S^\tau)$ is itself concentrated in
degree $k$ and thus $H^*_C(D^\tau,S^\tau)=H^k(D^\tau,S^\tau)[u]=\bbz_2[u]$.
The isomorphism $\kappa\co  H^{2k}(D,S)\to H^k(D^\tau,S^\tau)$ is obvious.
As $(D,S)$ is \ef, the restriction homomorphism
$r\co H^{*}_C(D,S)\to H^*_C(D^\tau,S^\tau)$ is injective.
This is a consequence of the localization theorem for singular cohomology,
which holds for smooth
actions on compact manifolds. Therefore,
if $a$ is the non-zero element of $H^{2k}(D,S)$, the equation
$r\sigma(a)=\kappa(a)u^k$ holds trivially. Hence, $(D,S)$
is a conjugation pair.
\end{Example}

\begin{Example}[Conjugation spheres]\label{spheres}
If $D$ is a conjugation cell of dimension $2k$ with boundary $S$,
the quotient space $\Sigma=D/S$ is a \csp\ homeomorphic to the sphere $S^{2k}$,
while $\Sigma^\tau$ is homeomorphic to $S^k$. For $a\in H^{2k}(\Sigma)$,
the conjugation equation $r\pcirc\sigma(a)=\kappa(a)u^k$ holds. 
We call such $\Sigma$ a {\it conjugation sphere}.
\end{Example}

\begin{Example}[Projective spaces]\label{excpn}
Let us consider the complex projective space $\bbc P^k$ with the involution
complex conjugation, having $\bbr P^k$ as real locus. 
One has $H^{2*}(\bbc P^k)=\bbz_2[a]/(a^{k+1})$ and 
and $H^*(\bbr P^k)=\bbz_2[b]/(b^{k+1})$.
The quotient space $\bbc P^k/\bbc P^{k-1}$ is a conjugation sphere.
Hence, in the following commutative diagram,
\begin{equation}\label{excpnEq1}
\xymatrix@C-3pt@M+2pt@R-4pt{%
H^{2k}_C(\bbc P^k,\bbc P^{k-1}) \ar[d]^(0.50){\approx}_(0.50){\rho_{\rm rel}}
\ar[r]^(0.60){\hat i}  &
H^{2k}_C(\bbc P^k) \ar@{>>}[d]^(0.50){\rho}  \\
%%%%%% ROW 2
H^{2k}(\bbc P^k,\bbc P^{k-1}) \ar[r]^(0.55){i}_(0.55){\approx}  &
H^{2k}(\bbc P^k)   \, ,
}
\end{equation}
the map $\rho_{\rm rel}$ is an isomorphism and hence
$\rho$ is surjective.
Setting $\sigma_{\rm rel}\co =\rho_{\rm rel}^\mun$, one gets a section
$\sigma$ of $\rho$ by $\sigma\co =\hat i\pcirc\sigma_{\rm rel}\pcirc i$.
The isomorphism
$$\kappa_{\rm
rel}\co H^{2*}(\bbc P^k,\bbc P^{k-1})\hfl{\approx}{} H^{*}(\bbr
P^k,\bbr P^{k-1})$$
is obvious and satisfies
$i^\tau\pcirc\kappa_{\rm rel}=\kappa\pcirc i$, where
$\kappa\co H^{2*}(\bbc P^k)\hfl{\approx}{} H^{*}(\bbr P^k)$ is the
unique ring isomorphism satisfying $\kappa(a)=b$.
Moreover, using $\hat\sigma_{\rm rel}$,
we have $H^*_C(\bbc P^k,\bbc P^{k-1})=H^{2k}(\bbc P^k,\bbc
P^{k-1})[u]$, and using the K\"unneth formula, $H^*(\bbr P^k,\bbr
P^{k-1})=H^{k}(\bbr P^k,\bbr P^{k-1})[u]$. Let $c\in H^{2k}(\bbc
P^k,\bbc P^{k-1})$ and $c'\in H^{k}(\bbr P^k,\bbr P^{k-1})$ be the
non-zero elements. As $\bbc P^k/\bbc P^{k-1}$ is a conjugation sphere,
the equation $r\pcirc\sigma_{\rm rel}(c)=c'u^k$ holds, giving
the formula $r\pcirc\sigma(a^k)=b^ku^k$ in $H^{2k}_C(\bbr P^k)$.

Now, if $k\leq n\leq\infty$, the restriction homomorphisms
$H^{2*}(\bbc P^n)\to H^{2*}(\bbc P^k)$,
$H^{2*}_C(\bbc P^n)\to H^{2*}_C(\bbc P^k)$,
$H^{*}(\bbr P^n)\to H^{*}(\bbr P^k)$ and
$H^{*}_C(\bbr P^n)\to H^{*}_C(\bbr P^k)$ are isomorphisms for $*\leq k$.
Therefore, the equation $r\pcirc\sigma(a^k)=b^ku^k$
holds in $H^{2k}_C(\bbr P^n)$  modulo
elements in the kernel of the restriction homomorphism
$H^{2k}_C(\bbr P^n)\to H^{2k}_C(\bbr P^k)$. This kernel
consists of terms of type $\lt{k}$. Therefore,
one has $r\pcirc\sigma(a^k)=b^ku^k+\lt{k}=\kappa(a^k)u^k+\lt{k}$
which shows that {\sl $\bbc P^n$ is a \csp\ for all $n\leq\infty$}.

We now show that the terms $\lt{k}$ in $H^{2k}_C(\bbr P^n)$ never
vanish when $n\geq 2k$. Let $\rho^\tau\co H^*_C(\bbr P^n)\to
H^*(\bbr P^n)$ and $r_0\co  H^*(\bbc P^n)\to H^*(\bbr P^n)$ be the
restriction homomorphisms. One has $\rho^\tau\pcirc r\pcirc\sigma=
r_0\pcirc\rho\pcirc\sigma$ and it is classical that $r_0(a)=b^2$
($a$ is the (${\rm mod\,} 2$) Euler class of the Hopf bundle
$\eta$ over $\bbc P^\infty$ and $b$ is the Euler class of the real Hopf bundle
$\eta^\tau$ over $\bbr P^\infty$; these bundles satisfy
$\eta_{|\bbr P^\infty}=\eta^\tau\oplus\eta^\tau$).
Therefore, $r(a)=bu+b^2$. Since $r\pcirc\sigma$ is a ring
homomorphism by Theorem~\ref{sigmaring}, one has
\begin{equation}\label{cfcp}
r\pcirc\sigma(a^k)= (bu + b^2)^k \, .
\end{equation}
Therefore, a term $b^{2k}$ is always present in the right hand
side of~\eqref{cfcp} when $n\geq 2k$. For instance,
$r\pcirc\sigma(a^2)=b^2u^2+b^4$, $r\pcirc\sigma(a^3)=b^3u^3+
b^4u^2 + b^5u + b^6$, and so on.
\end{Example}

We finish this section with two related results.

\begin{Lemma}[Injectivity lemma]\label{injthm}
Let $(X,Y)$ be a \cpa. Then the restriction homomorphism
$r\co H^*_C(X,Y)\to H^*_C(X^\tau,Y^\tau)$ is injective.
\end{Lemma}

\begin{proof}
Suppose that $r$ is not injective.
Let $0\neq x=\sigma(y)u^k\! + \lt{k}{\in} H^{2n+k}_C(X,Y)$
be an element in $\ker r$. The conjugation equation guarantees
that $k\neq 0$. We may assume that $k$ is minimal.
By the conjugation equation again,
we have $0=r(x)=\kappa(y)u^{n+k}+\lt{n+k}$. Since $\kappa$ is an isomorphism,
we get $y=0$, which is a contradiction. 
\end{proof}

\begin{Lemma}\label{locthm}
Let $(X,Y)$ be a \cpa. Assume that $H^{2*}(X,Y)=0$ for $*>m_0$.
Then the localization theorem holds. That is, the restriction homomorphism
$r\co H^*_C(X,Y)\to H^*_C(X^\tau,Y^\tau)$ becomes an isomorphism
after inverting $u$.
\end{Lemma}

\begin{proof} 
By Lemma~\ref{injthm}, it suffices to show that
$$
H^*(X^\tau,Y^\tau)=H^*(X^\tau,Y^\tau)\otimes 1 \subset H^*_C(X^\tau,Y^\tau)
$$
is in the image of $r$ localized. We show this by downward induction
on the degree of an element in $H^*(X^\tau,Y^\tau)$. The statement is
obvious for $*>m_0$. Since $r\pcirc\sigma(x)=\kappa(x)u^k+\lt{k}$ 
for $x\in H^{2k}(X,Y)$, the induction step follows 
(by induction hypothesis, $\lt{k}$ is in the image of $r$ localized). 
\end{proof}

\begin{Remark}\label{locinj}
In classical equivariant cohomology theory, the injectivity lem\-ma is often
deduced from the localization theorem. But, as seen in Example~\ref{excpn},
$\bbc P^\infty$ with the complex conjugation is a \csp, and therefore
satisfies the injectivity lemma. However,
$r_{\rm loc} \co  H^*_C(\bbc P^\infty)[u^\mun]\to H^*_C(\bbr P^\infty)[u^\mun]$
is not surjective. Indeed, $H^*_C(\bbc P^\infty)[u^\mun]=\bbz_2[a,u,u^\mun]$,
$H^*_C(\bbr P^\infty)[u^\mun]=\bbz_2[b,u,u^\mun]$ and $r_{\rm loc}(a)=bu+b^2$
by Example~\ref{excpn}. Therefore, $r_{\rm loc}$ composed with the epimorphism
$\bbz_2[b,u,u^\mun]\to\bbz_2$ sending $b$ and $u$ to $1$ is the zero map.
\end{Remark}

%%%%%%%%%%%%%%%%%%%%%%%%%%%%%%%%%%%%%%%%%%%%%%%%%%%%%%%%
\subsection{Equivariant maps between \csp s}\label{Seqmapsalc}

The purpose of this section is to show the naturality of $H^*$-frames.
Let $(X,X_0)$ and $(Y,Y_0)$ be two \cpa s.
Choose \hfra s $(\kappa_{\scr X},\sigma_{\scr X})$
and  $(\kappa_{\scr Y},\sigma_{\scr Y})$
for $(X,X_0)$ and $(Y,Y_0)$ respectively.
Let $f\co (Y,Y_0)\to (X,X_0)$ be a $C$-equivariant
map of pairs.
We denote by $f^\tau\co  (Y^\tau,Y_0^\tau)\to (X^\tau,X^\tau_0)$
the restriction of $f$ to $(Y^\tau,Y_0^\tau)$
and use the functorial notations : $H^*f$, $H^*_Cf$, and so forth.

\begin{Proposition}\label{alchf}
The \csp\ structure of a \csp\ is natural, i.e., one has
\begin{equation}\label{alchf-eqA}
H^*_Cf\pcirc\sigma_{\scr X} = \sigma_{\scr Y}\pcirc H^*f
\end{equation}
and
\begin{equation}\label{alchf-eqB}
H^*f^\tau\pcirc\kappa_{\scr X} = \kappa_{\scr Y}\pcirc H^*f \ .
\end{equation}
\end{Proposition}

\begin{proof}
Let $\rho_{\scr X}\co H^{*}_C(X,X_0)\to H^{*}(X,X_0)$
and $\rho_{\scr Y}\co H^{*}_C(Y,Y_0)\to H^{*}(Y,Y_0)$
denote the restriction homomorphisms.
Let $a\in H^{2k}(X,X_0)$.
As $H^*f\pcirc\rho_{\scr X}=
\rho_{\scr Y}\pcirc H^*_Cf$, one has
\begin{equation}\label{alchf-eq1}
\rho_{\scr Y}\pcirc H^*_Cf\pcirc\sigma_{\scr  X}(a) =
H^*f\pcirc\rho_{\scr X}\pcirc\sigma_{\scr  X}(a)  =
H^*f(a) =
\rho_{\scr Y}\pcirc\sigma_{\scr  Y}\pcirc H^*f(a) .
\end{equation}
This implies that Equation~\eqref{alchf-eqA} holds modulo the ideal
$(u)$. As $H^*(X,X_0)$ is concentrated in even degrees, this means that
\begin{equation}\label{alchf-eq2}
H^*_Cf\pcirc\sigma_{\scr X}(a) = \sigma_{\scr  Y}\pcirc H^*f(a) +
\sigma_{\scr  Y}(d_{2k-2})u^2 + \cdots + \sigma_{\scr  Y}(d_{0})u^{2k} \ ,
\end{equation}
where $d_i\in H^{i}(Y,Y_0)$.
Now, by Lemma \ref{hfpri},
\begin{equation}\label{alchf-eq30}
H^*_Cf^\tau\pcirc r_X\pcirc\sigma_{\scr X}(a) = H^*_Cf^\tau(\kappa_{\scr X}(a)\,u^k+\lt{k}) =
H^*f^\tau\pcirc\kappa_{\scr X}(a)\,u^k  + \lt{k}.
\end{equation}
On the other hand, by equation \eqref{alchf-eq2}
\begin{equation}\label{alchf-eq40}
r_Y\pcirc H^*_Cf(a) = \sigma_{\scr Y}(d_0) u^{2k} + \lt{2k}.
\end{equation}
But $r_Y\pcirc H^*_Cf= H^*_Cf^\tau\pcirc r_X$.
Comparing then Equation~\eqref{alchf-eq40} with
Equation~\eqref{alchf-eq30}, we deduce that
$d_0=0$, since $\kappa_{\scr Y}$ is an injective. Then
\begin{equation}\label{hf-eq50}
r_Y\pcirc H^*_Cf(a) =
\kappa_{\scr Y}(d_2)\,u^{2k-2}
+ \lt{2k-2}.
\end{equation}
Again, we deduce that $d_2=0$.
Continuing this process, we finally get Equation~\eqref{alchf-eqA}
(as in the proof of Theorem~\eqref{sigmaring}).

As for Equation~\eqref{alchf-eqB}, by Lemma~\ref{hfpri},
\begin{equation}\label{alceqmapcoh-eqa}
H^*_Cf^\tau\pcirc r_X\pcirc\sigma_{\scr X}(a) =
H^*_Cf^\tau(\kappa_{\scr X}(a)u^k + \lt{k}) =
H^*f^\tau\pcirc\kappa_{\scr X}(a)u^k + \lt{k} \ .
\end{equation}
On the other hand, using Equation~\eqref{alchf-eqA},
\begin{equation}\label{alceqmapcoh-eqb}
r_Y\pcirc H^{*}_Cf\pcirc\sigma_{\scr X}(a) =
r_Y\pcirc\sigma_{\scr Y}\pcirc H^{*}f(a)
=\kappa_{\scr Y}\pcirc H^{*}f(a)u^k + \lt{k} \ .
\end{equation}
Comparing Equation~\eqref{alceqmapcoh-eqa} with \eqref{alceqmapcoh-eqb} gives
Equation~\eqref{alchf-eqB}. 
\end{proof}

As a corollary of \proref{alchf}, we get the uniqueness
of the \csp\ structure for a \csp .

\begin{Corollary}\label{alcstrunique}
Let $(\kappa,\sigma)$ and $(\kappa',\sigma')$ be two
\hfra s for an involution $\tau$ on $(X,X_0)$ Then
$(\kappa,\sigma)=(\kappa',\sigma')$
\end{Corollary}

\begin{proof} 
If two \hfra s
$(\kappa_{\scr X},\sigma_{\scr X})$ and $(\kappa_{\scr X}',\sigma_{\scr X}')$
are given on $(X,X_0)$, \proref{alchf} with $f={\rm id}_X$ proves
that $\kappa_{\scr X}=\kappa_{\scr X}'$ and
$\sigma_{\scr X}=\sigma_{\scr X}'$. 
\end{proof}

By the Leray-Hirsch Theorem, the section 
$\sigma_{\scr  X}\co  H^{*}(X,X_0)\to H^{*}_C(X,X_0)$  induces a map
$\hat\sigma_{\scr  X}\co H^{*}(X,X_0)[u] \stackrel{\approx}\to H^{*}_C(X,X_0)$
which is an isomorphism of $R$-algebras by Corollary \ref{sigmaRalgebrabis}. 
We define 
$\hat\sigma_{\scr  Y}\co H^{*}(Y,Y_0)[u] \stackrel{\approx}\to H^{*}_C(Y,Y_0)$
accordingly. Proposition \ref{alchf} shows that these $R$-algebras
isomorphisms are natural and gives the following analogue of
Lemma \ref{hfpri}.
\begin{Corollary}\label{alcsigmaalgnat}
For any $C$-equivariant map $f\co Y\to X$ between \csp s,
the diagram
$$\xymatrix@C-3pt@M+2pt@R-4pt{%
H^*(X,X_0)[u] \ar[r]^(0.55){\hat\sigma_{\scr  X}}_(0.55){\approx} \ar[d]_{H^*f[u]}
& H^*_C(X,X_0)\ar[d]^{H^*_Cf} \\
H^*(Y,Y_0)[u] \ar[r]^(0.55){\hat\sigma_{\scr  Y}}_(0.55){\approx} & H^*_C(Y,Y_0)
}$$
is commutative, where $H^*f[u]$
is the polynomial extension of $H^*f$.
 \qed
\end{Corollary}

Finally, \proref{alchf} and Corollary~\ref{alcsigmaalgnat}
give the naturality of the algebra isomorphism $\hat\kappa$
of Equation~\eqref{defhatkappa}.

\begin{Proposition}\label{eqmapcoh}
For any $C$-equivariant map $f\co Y\to X$ between \csp s,
the diagram
$$
\begin{array}{cccccc}
H^{2*}_C(X,X_0)  & \hfl{H^{2*}_Cf}{} & H^{2*}_C(Y,Y_0) \\
\downarrow \hat\kappa_{\scr  X} && \downarrow \hat\kappa_{\scr  Y} \\
H^*_C(X^\tau,X_0^\tau) & \hfl{H^*_Cf^\tau}{} & H^*_C(Y^\tau,Y_0^\tau)
\end{array}
$$
is commutative. \qed
\end{Proposition}

%%%%%%%%%%%%%%%%%%%%%%%%%%%%%%%%%%%%%%%%%%%%%%%%%%%%%%%%
\section{Extension properties}\label{SEXP}

\subsection{Triples}\label{SStriples}

\begin{Proposition}\label{EXPlac}
Let $X$ be a space with an involution $\tau$ and let $Z\subset Y$ be
$\tau$-invariant subspaces of $X$. Suppose that
$(X,Y)$ and $(Y,Z)$ are \cpa s.
Then $(X,Z)$ is a \cpa.
\end{Proposition}

\begin{proof}
The subscript ``$X,Y$'' is used for the 
relevant homomorphism for the 
pair $(X,Y)$, like $\kappa_{\scr X,Y}$, $r_{\scr X,Y}$, etc.
In order to simplify the notation, we use the subscripts 
``$X$'' or ``$Y$'' for the pairs $(X,Z)$ and $(Y,Z)$, 
as if $Z$ were empty.
Thus, we must construct a \hfra\ $(\kappa_{\scr X},\sigma_{\scr X})$
for the pair $(X,Z)$, using those $(\kappa_{\scr Y},\sigma_{\scr Y})$
and $(\kappa_{\scr X,Y},\sigma_{\scr X,Y})$ for the \cpa s
$(Y,Z)$ and $(X,Y)$.
 
We first prove that the restriction homomorphisms
$\hat j\co H^*_C(X,Z)\to H^*_C(Y,Z)$ and
$j^\tau\co H^*(X^\tau,Z^\tau)\to H^*(Y^\tau,Z^\tau)$ are surjective.
Let us consider the following commutative diagram
\begin{equation}\label{EXPlac-eq501}
\xymatrix@C-3pt@M+2pt@R-4pt{%
H^*_C(Y,Z) \ar@{>->}[d]_(0.50){r_Y}
\ar[r]^(0.50){\delta_{\scr C}}  &
H^{*+1}_C(X,Y)  \ar@{>->}[d]^(0.50){r_{\scr X}} \\
%%%%%% ROW 2
H^*_C(Y^\tau,Z^\tau) \ar[r]^(0.50){\delta^\tau_{\scr C}}  &
H^{*+1}_C(X^\tau,Y^\tau)
}
\end{equation}
in which $\delta_{\scr C}$ and $\delta^\tau_{\scr C}$ are the connecting homomorphisms
for the long exact sequences in equivariant cohomology of the triples
$(X,Y,Z)$ and $(X^\tau,Y^\tau,Z^\tau)$ respectively.
The vertical restriction homomorphisms are injective by Lemma~\ref{injthm}.
Clearly, $\delta_{\scr C}=0$ if and only if $\hat j$ is surjective.
As $\delta^\tau_{\scr C}$ is the polynomial extension of $\delta^\tau$,
one also has
$$\delta^\tau_{\scr C}=0 \ \Longleftrightarrow \ \delta^\tau=0
\ \Longleftrightarrow \ \hbox{$j^\tau$ is surjective}\, .
$$
As $(X,Y)$ is an even cohomology pair,
for $y\in H^{2k}(Y,Z)$, one can write
\begin{equation}\label{EXPlac-eq510}
\delta_{\scr C}\pcirc\sigma_{\scr Y}(y)=
\sum_{i=0}^k \sigma_{\scr X}(x_{2k-2i})u^{2i+1} \, ,
\end{equation}
with $x_{2k-2i}\in H^{2k-2i}(X,Y)$.
Using that $r_{\scr Y}\pcirc \sigma_{\scr Y}(y) = \kappa_{\scr Y}(y)u^k+\lt{k}$,
the commutativity of Diagram~\eqref{EXPlac-eq501} and that
$\delta^\tau_{\scr C}=\delta^\tau[u]$, we get
\begin{equation}\label{EXPlac-eq520}
\begin{array}{rcl}
\delta^\tau\pcirc\kappa_{\scr Y}(y)u^k +\lt{k}  &=&
\delta^\tau_{\scr C}(\kappa_{\scr Y}(y)u^k+\lt{k}) \\[2mm] &=&
r_{\scr X}\big(\displaystyle\sum_{i=0}^k \sigma_{\scr X}(x_{2k-2i})u^{2i+1}\big)
\\[2mm] &=&
\displaystyle\sum_{i=0}^k\big(\kappa_{\scr X}(x_{2k-2i})u^{k+i+1} + \lt{k+i+1}\big) \ .
\end{array}
\end{equation}
As in the proof of Theorem~\ref{sigmaring},
we compare the coefficients of powers of $u$,
in both sides of Equation~\ref{EXPlac-eq520}.
Starting with $u^{2k+1}$ and going downwards, 
we get inductively that $\kappa_{\scr X}(x_{2k-2i})=0$ for $i=0,\dots,k$.
Hence $x_{2k-2i}=0$ for $i=0,\dots,k$ and the right side of
Equation~\ref{EXPlac-eq520} vanishes for all $y\in H^{2k}(Y,Z)$.
As $\kappa_{\scr Y}$ is bijective, we deduce that $\delta^\tau=0$ and
$\delta^\tau_{\scr C}=0$. As $r_{\scr X}$ is injective by Lemma~\ref{injthm},
the commutativity of Diagram~\eqref{EXPlac-eq501} implies that $\delta_{\scr C}=0$.
We have thus proven that the restriction homomorphisms
$\hat j\co H^*_C(X,Z)\to H^*_C(Y,Z)$ and
$j^\tau\co H^*(X^\tau,Z^\tau)\to H^*(Y^\tau,Z^\tau)$ are surjective.

As $\rho_{\scr Y}$ is onto, the cohomology exact sequence of $(X,Y,Z)$
decomposes into short exact sequences and one has
the following commutative diagram:
\begin{equation}\label{EXPlac-eq10}
\xymatrix@C-3pt@M+2pt@R-4pt{%
0 \ar[r] &  H^{2*}_C(X,Y) \ar[r]^{\hat i} \ar@{>>}[d]^{\rho_{\scr X,Y}} &
H^{2*}_C(X,Z) \ar[r]^{\hat j} \ar[d]^{\rho_{\scr X}} &
H^{2*}_C(Y,Z) \ar@{>>}[d]^{\rho_{\scr Y}} \ar[r]\ar@/^/[l]^{\hat\mu}
\ar[r] & 0
\\
0 \ar[r] &
H^{2*}(X,Y) \ar[r]^i    &
H^{2*}(X,Z) \ar[r]^j  &
H^{2*}(Y,Z)\ar[r]\ar@/^/[l]^{\mu} \ar[r] & 0
}
\end{equation}
Sections $\hat\mu$ and $\mu$ can be constructed as follows.
Let $\calb$ be a basis of the $\bbz_2$-vector space
$H^{2*}(Y,Z)$. The set $\sigma_{\scr Y}(B)$ is a
$R^{ev}$-module basis for $H^{2*}_C(Y,Z)$.
For each $b\in\calb$, choose
$\tilde b\in H^{2*}_C(X,Z)$
such that $\hat j(\tilde b)=\sigma_{\scr  Y}(b)$. The correspondence $\sigma_{\scr  Y}(b)\to \tilde b$
induces a section $\hat\mu\co  H^{2*}_C(Y,Z)\to H^{2*}_C(X,Z)$ of
$\hat j$. One has $j\pcirc\rho_{\scr X}\pcirc\hat\mu\pcirc\sigma_{\scr  Y}(b)=b$;
therefore $\mu\co =\rho_{\scr X}\pcirc\hat\mu\pcirc\sigma_{\scr Y}$
is an additive section of the epimorphism $j$.

Using that additively, $H^{2*}(X,Z)=i(H^{2*}(X,Y))\oplus\mu(H^{2*}(Y,Z))$,
one defines $\sigma_{\scr X}^0\co H^{2*}(X,Z)\to H^{2*}_C(X,Z)$ by:
\begin{equation}\label{EXPlac-eq20}
\left\{
\begin{array}{rcll}
\sigma_{\scr X}^0(i(a)) &\ccoeg &
\hat i\pcirc\sigma_{\scr X,Y}(a) & \hbox{ for all }
a\in H^{2*}(X,Y)\\[2mm]
\sigma_{\scr X}^0(\mu(b))&\ccoeg &
\hat\mu\pcirc\sigma_{\scr Y}(b)& \hbox{ for all }
b\in H^{2*}(Y,Z)
\end{array}\right.
\end{equation}
The map $\sigma_{\scr X}^0$ is an additive section of $\rho_{\scr X}$
and the following diagram is commutative:
\begin{equation}\label{EXPlac-eq30}
\xymatrix@C-3pt@M+2pt@R-4pt{%
0 \ar[r] &  H^{2*}(X,Y) \ar[r]^{i} \ar@{>->}[d]^{\sigma_{\scr X,Y}} &
H^{2*}(X,Z) \ar[r]^{j} \ar@{>->}[d]^{\sigma_{\scr X}^0} &
H^{2*}(Y,Z) \ar@{>->}[d]^{\sigma_{\scr Y}}
\ar[r] & 0
\\
0 \ar[r] &
H^{2*}_C(X,Y) \ar[r]^{\hat i}    &
H^{2*}_C(X,Z) \ar[r]^{\hat j}  &
H^{2*}_C(Y,Z)\ar[r] \ar[r] & 0
}
\end{equation}
We define an additive map 
$\kappa_{\scr X}^0\co H^{2*}(X,Z)\to H^*(X^\tau,Z^\tau)$ by
\begin{equation}\label{EXPlac-eq44}
\left\{
\begin{array}{rcll}
\kappa_{\scr X}^0(i(a)) &\ccoeg &
i^\tau\pcirc\kappa_{\scr X,Y}(a) & \hbox{ for all }
a\in H^{2*}(X,Y)\\[2mm]
\kappa_{\scr X}^0(\mu(b))&\ccoeg &
\mu^\tau\pcirc\kappa_{\scr Y}(b)& \hbox{ for all }
b\in H^{2*}(Y,Z) \, ,
\end{array}\right.
\end{equation}
where $\mu^\tau\co H^{*}(Y^\tau,Z^\tau)\to H^{*}(X^\tau,Z^\tau)$ is any additive
section of $j^\tau$.
The following diagram is then commutative:
\begin{equation}\label{EXPlac-eq45}
\xymatrix@C-3pt@M+2pt@R-4pt{%
0 \ar[r] &  H^{2*}(X,Y) \ar[r]^{i} \ar[d]_{\kappa_{\scr X,Y}}^{\approx} &
H^{2*}(X,Z) \ar[r]^{j} \ar[d]^{\kappa_{\scr X}^0} &
H^{2*}(Y,Z) \ar[d]^{\kappa_{\scr Y}}_{\approx}
\ar[r] & 0
\\
0 \ar[r] &
H^{*}(X^\tau,Y^\tau) \ar[r]^{i^\tau}    &
H^{*}(X^\tau,Z^\tau) \ar[r]^{j^\tau}  &
H^{*}(Y^\tau,Z^\tau)\ar[r] \ar[r] & 0
}\end{equation}
By construction, the equality
$r_X\pcirc\sigma_{\scr X}^0(i(a))=\kappa_{\scr X}^0(i(a))\,u^k + \lt{k}$
holds for all $a\in H^{2k}(X,Y)$ and all $k$.
On the other hand, for $b\in H^{2k}(Y,Z)$, we only have that
$j^\tau\pcirc r_X\pcirc\sigma_{\scr X}^0(\mu(b))=
j^\tau\big(\kappa_{\scr X}^0(\mu(b))\,u^k + \lt{k})$, which implies that
\begin{equation}\label{EXPlac-eq50}
r_X\pcirc\sigma_{\scr X}^0(\mu(b))=
\hat i^\tau(D) + \kappa_{\scr X}^0(\mu(b))\,u^k + \lt{k}
\end{equation}
for some $D\in H^{*}_C(Y^\tau,Z^\tau)$. As $\hat i^\tau$ is injective
(since $\delta^\tau_C$ is onto), the element $D$ in
Equation~\eqref{EXPlac-eq50} is unique if chosen free of terms
$\lt{k}$. Such a $D$ is of the form
\begin{equation}\label{EXPlac-eq51}
D=\kappa_{\scr X,Y}(d_{2k})u^k + \sum_{s=1}^k \kappa_{\scr X,Y}(d_{2(k-s)})u^{k+s} \, ,
\end{equation}
where $d_i\in H^{i}(X,Y)$. Define
$\sigma_{\scr X}\co H^{2*}(X,Z)\!\to\! H^{2*}_C(X,Z)$
and
$\kappa_{\scr X}\co H^{2*}(X,Z)$ $\to H^*(X^\tau,Z^\tau)$ by
$$\sigma_{\scr X}(i(a)) \ccoeg  \sigma_{\scr X}^0(i(a))
\ \hbox{ and }\
\kappa_{\scr X}(i(a)) = \kappa_{\scr X}^0(i(a)) \
\hbox{ for all } a \in H^{2*}(X,Y),$$
and, for $b\in H^{2*}(Y,Z)$, by
$$
\left\{
\begin{array}{rcll}
\sigma_{\scr X}(\mu(b))&\ccoeg & \sigma_{\scr X}^0(\mu(b))
+\sum_{s=1}^k i(d_{2(k-s)})u^{2s}
\\[2mm]
\kappa_{\scr X}(\mu(b))&\ccoeg & \kappa_{\scr X}^0(\mu(b))
+i^\tau(d_{2k})
\end{array}\right. \ .$$
We may check that 
$r_X\pcirc\sigma_{\scr X}(c)=\kappa_{\scr X}(c)\,u^k + \lt{k}$ 
for all $c\in H^{2k}(X,Z)$.
As $\sigma_{\scr X}(c)-\sigma_{\scr X}^0(c)\in H^{*}(X,Z)\cdot u
=\ker\rho_{\scr X}$, the homomorphism $\sigma_{\scr X}$ is a section of $\rho_{\scr X}$. Diagram~\eqref{EXPlac-eq45} still commutes
with $\kappa_{\scr X}$ instead of $\kappa_{\scr X}^0$. As
$\kappa_{\scr X,Y}$ and $\kappa_{\scr Y}$ are bijective,
$\kappa_{\scr X}$ is bijective by the five-lemma. 
\end{proof}

\begin{Proposition}\label{EXPlac2}
Let $X$ be a space with an involution $\tau$ and let $Z\subset Y$ be
$\tau$-invariant subsets of $X$. Suppose that
\renewcommand{\labelenumi}{\rm(\roman{enumi})}
\begin{enumerate}
\item $(X,Z)$ and $(X,Y)$ are \cpa s.
\item the restriction homomorphisms $i\co H^*(X,Y)\to H^*(X,Z)$ 
is injective.
\end{enumerate}
Then $(Y,Z)$ is a \cpa.
\end{Proposition}

\begin{Remark}\label{EXPlacRem2} \rm Assuming condition (i), 
condition (ii) is necessary for $(Y,Z)$ to be a \cpa, 
since the three pairs will then have cohomology 
only in even degrees.
\end{Remark}

\begin{proof}[Proof of Proposition \ref{EXPlac2}]
We have the following commutative diagram
\begin{equation}\label{EXPlac2-eq10}
\xymatrix@C-3pt@M+2pt@R-4pt{%
&  H^{2*}_C(X,Y) \ar[r]^{\hat i} \ar[d]^{\rho_{\scr X,Y}} &
H^{2*}_C(X,Z) \ar[r]^{\hat j} \ar@{>>}[d]^{\rho_{\scr X}} &
H^{2*}_C(Y,Z) \ar[d]^{\rho_{\scr Y}}
\\
0 \ar[r] &
H^{2*}(X,Y) \ar[r]^i   \ar@/^/[u]^{\sigma_{\scr X,Y}} &
H^{2*}(X,Z) \ar[r]^j \ar@/^/[u]^{\sigma_{\scr X}} &
H^{2*}(Y,Z)\ar[r]\ar@/^/[l]^{\mu} \ar[r] & 0
}
\end{equation}
where $\mu$ is an additive section of $j$.
Define an additive section $\sigma_{\scr Y}$ of $\rho_{\scr Y}$
by $\sigma_{\scr Y}\ccoeg \hat j\pcirc\sigma_{\scr X}\pcirc\mu$ and
$\kappa_{\scr Y}\co H^{2*}(Y,Z)\to H^{*}(Y^\tau,Z^\tau)$ by
$\kappa_{\scr Y}\ccoeg \hat j^\tau\pcirc\kappa_{\scr X}\pcirc\mu$.
This guarantees that
$r_{\scr Y}\sigma_{\scr Y}(a)=\kappa_{\scr Y}(a)u^k+\lt{k}$ for all
$a\in H^{2k}(Y,Z)$. It then just remains to prove that $\kappa_{\scr Y}$
is bijective. 

As $i$ is injective,
the equation $i^\tau\pcirc\kappa_{\scr X,Y}=\kappa_{\scr X}\pcirc i$,
guaranteed by \proref{alchf}, implies that $i^\tau$ is injective.
The same equation implies that 
$j^\tau\pcirc\kappa_{\scr X}=\kappa_{\scr Y}\pcirc j$,
since $j^\tau\pcirc\kappa_{\scr X}\pcirc i=0$. Therefore,
one has a commutative diagram
\begin{equation}\label{EXPlac2-eq45}
\xymatrix@C-3pt@M+2pt@R-4pt{%
0 \ar[r] &  H^{2*}(X,Y) \ar[r]^{i} \ar[d]^{\kappa_{\scr X,Y}}_{\approx} &
H^{2*}(X,Z) \ar[r]^{j} \ar[d]^{\kappa_{\scr X}}_{\approx} &
H^{2*}(Y,Z) \ar[d]^{\kappa_{\scr Y}}
\ar[r] & 0
\\
0 \ar[r] &
H^{*}(X^\tau,Y^\tau) \ar[r]^{i^\tau}    &
H^{*}(X^\tau,Z^\tau) \ar[r]^{j^\tau}  &
H^{*}(Y^\tau,Z^\tau)\ar[r] \ar[r] & 0  \, 
}
\end{equation}
which shows that $\kappa_Y$ is bijective.
\end{proof}

The same kind of argument will prove \proref{EXPlac3} below. As this
proposition is not used elsewhere in this paper, we leave the proof
to the reader.

\begin{Proposition}\label{EXPlac3}
Let $X$ be a space with an involution $\tau$ and let $Z\subset Y$ be
$\tau$-invariant subsets of $X$. Suppose that
\renewcommand{\labelenumi}{\rm(\roman{enumi})}
\begin{enumerate}
\item $(X,Z)$ and $(Y,Z)$ are \cpa s.
\item the restriction homomorphisms 
$j\co H^*(X,Z)\to H^*(Y,Z)$ is surjective.
\end{enumerate}
Then $(X,Y)$ is a \cpa. \qed
\end{Proposition}

%%%%%%%%%%%%%%%%%%%%%%%%%%%%%%%%%%%%%%%%%%%%%%%%%
\subsection{Products}\label{Sprod}

\begin{Proposition}\label{Pprod}
Let $(X,X_0)$ and $(Y,Y_0)$ be \cpa s.
Suppose that $H^q(X,X_0)$ is finite dimensional for each $q$. Assume that
$\{X\times Y_0,X_0\times Y\}$ is an excisive couple in $X\times Y$ and that
$\{X^\tau\times Y_0^\tau,X_0^\tau\times Y^\tau\}$ is an excisive couple in 
$X^\tau\times Y^\tau$. Then, the product pair 
$(X\times Y,(X_0\times Y)\cup (X\times Y_0))$ is a \cpa. 
\end{Proposition}

\begin{proof} 
To simplify the notations, we give the proof when $X_0=Y_0=\emptyset$;
the general case is identical.
By of our hypotheses, the two projections $X\times Y\to X$ 
and $X\times Y\to Y$ give rise
to the K\"unneth isomorphism
$$K\co H^*(X)\otimes H^*(Y)\hfl{\approx}{} H^*(X\times Y) \, .$$
The same holds for the fixed point sets, producing
$$K^\tau\co H^*(X^\tau)\otimes H^*(Y^\tau)\hfl{\approx}{}
H^*(X^\tau\times Y^\tau) = H^*((X\times Y)^\tau) \, .$$

The Borel construction applied to the projections gives rise
to maps $(X\times Y)_C\to X_C$ and  $(X\times Y)_C\to Y_C$. This
produces a ring homomorphism
$$K_C\co H^*(X_C)\otimes H^*(Y_C)\hfl{}{} H^*((X\times Y)_C) \, .$$
We now want to define
$\kappa_{\scr X\times Y}$ and $\sigma_{\scr X\times Y}$.  We
set $\kappa_{\scr X\times Y} \ccoeg  K^\tau \pcirc (\kappa_{\scr
X}\otimes\kappa_{\scr Y}) \pcirc K^\mun$.  Then,
$\kappa_{\scr X\times Y}$ is an isomorphism and one has
the following commutative diagram:
$$
\xymatrix@C-3pt@M+2pt@R-4pt{%
H^{2*}(X)\otimes H^{2*}(Y)
\ar[d]_(0.50){\kappa_{\scr X}\otimes\kappa_{\scr Y}}^{\approx}
\ar[r]^(0.55){K}_(0.55){\approx}  &
H^{2*}(X\times Y) \ar[d]^(0.50){\kappa_{\scr X\times Y}}  \\
%%%%%% ROW 2
H^{*}(X^\tau)\otimes H^{*}(Y^\tau)
\ar[r]^(0.55){K^\tau}_(0.55){\approx}  & H^{*}((X\times Y)^\tau)
}
$$
Now, setting $\sigma_{\scr X\times Y}\ccoeg  K_C \pcirc (\sigma_{\scr
X}\otimes\sigma_{\scr Y}) \pcirc K^\mun$, we have:
$$
\xymatrix@C-0pt@M+2pt@R-4pt{%
H^{2*}_C(X)\otimes H^{2*}_C(Y) \ar[d]_(0.50){\rho_{\scr X}\otimes\rho_{\scr Y} }
\ar[r]^(0.55){K_C}  &
H^{2*}_C(X\times Y)
\ar[d]_(0.50){\rho_{\scr X\times Y}}  \\
%%%%%% ROW 2
H^{2*}(X)\otimes H^{2*}(Y) \ar@/_/[u]_{\sigma_{\scr
X}\otimes\sigma_{\scr Y}} \ar[r]^(0.55){K}_(0.55){\approx}  &
H^{2*}(X\times X) \ar@/_/[u]_{\sigma_{\scr X\times Y}}
}
$$
With these definitions, one verifies the conjugation 
equation by direct computation. 
\end{proof}

%%%%%%%%%%%%%%%%%%%%%%%%%%%%%%%%%%%%%%%%%%%%%%%%%%
\subsection{Direct limits}\label{colim}

\begin{Proposition}\label{Pcolim}
Let $(X_i,f_{ij})$ be a directed system of \csp s and $\tau$-equivariant
inclusions, indexed by a direct set $\cali$. Suppose that each space $X_i$
is $T_1$. Then $X=\displaystyle\lim_\to X_i$ is a \csp.
\end{Proposition}

\begin{proof} 
As the maps $f_{ij}$ are inclusion between and each $X_i$ is $T_1$,
the image of a compact set $K$ under a continuous
map to $X$ is contained in some $X_i$ (otherwise $K$ would
contain an infinite closed discrete subspace). Therefore,
$H_*(X)=\lim_\to H_*(X_i)$
(singular homology with $\bbz_2$ as coefficients).
Then
\begin{equation}\label{Pcolim-eq1}
\begin{array}{rcl}
H^*(X)={\rm Hom}(H_*(X);\bbz_2) &=&
{\rm Hom}\big(\displaystyle\lim_\to H_*(X_i);\bbz_2)\\ &=&
\displaystyle\lim_\leftarrow {\rm Hom}(H_*(X_i);\bbz_2)\\ &=&
\displaystyle\lim_\leftarrow H^*(X_i)\ .
\end{array}
\end{equation}
One has $X^\tau=\lim_\to X_i^\tau$ and $X_C=\lim_\to (X_i){}_C$ and, as
in~\eqref{Pcolim-eq1}, one has $H^*(X^\tau)=\lim_\leftarrow H^*(X_i^\tau)$
and $H^*_C(X)=\lim_\leftarrow H^*_C(X_i)$. By \proref{alchf},
the isomorphisms $\kappa_i \co  H^{2*}(X_i)\to H^*(X_i^\tau)$ is an isomorphism
of inverse systems; we can thus define $\kappa=\lim_\leftarrow \kappa_i$,
and $\kappa$ is an isomorphism. The same can be done for
$\sigma\co H^{2*}(X)\to H^{2*}_C(X)$, defined, using \proref{alchf},
as the inverse limit of $\sigma_i\co H^{2*}(X_i)\to H^{2*}_C(X_i)$,
and $\sigma$ is a section of $\rho\co H^{2*}_C(X)\to H^{2*}(X)$. The
conjugation equation for $(\sigma,\kappa)$ comes directly from that for $(\sigma_i,\kappa_i)$. 
\end{proof}

%%%%%%%%%%%%%%%%%%%%%%%%%%%%%%%%%%%%%%%%%%%%%%%%%%
\subsection{Equivariant connected sums}\label{ccsum}
Let $M$ be a smooth oriented closed manifold of dimension $2k$
together with a smooth involution $\tau$ that is a conjugation.
Then $M^\tau$ is a non-empty closed submanifold of $M$ of
dimension $k$. Pick a point $p\in M^\tau$. There is a
$\tau$-invariant disk $\Delta$ of dimension $2k$ in $M$ around $p$
on which $\tau$ is conjugate to a linear action: there is a
diffeomorphism $h\co \bbd(\bbr^k\times\bbr^k)\to \Delta$ preserving
the orientation such that $\tau\pcirc h = h\pcirc\tau_0$, where
$\tau_0(x,y)=(x,-y)$.

Let $(M_i,\tau_i)$, $i=1,2$, be two smooth conjugation spaces, as
above. Choosing conjugation cells 
(see Example~\ref{disk})
$h_i\co \bbd(\bbr^k\times\bbr^k)\to
\Delta_i$ as above, one can form the connected sum
$$M\ccoeg M_1\sharp M_2 = (M_1\setminus {\rm int}\Delta_1)\cup_{h_2\pcirc h_1^\mun}
(M_2\setminus {\rm int}\Delta_2)$$
which inherits
an involution $\tau$. We do not know whether the
equivariant diffeomorphism
type of $M_1\sharp M_2$ depends on the choice of
the diffeomorphism $h_i$, which is unique only up
to pre-composition by elements of $S(O(k)\times O(k))$.

\begin{Proposition}\label{connsum}
$M_1\sharp M_2$ is a conjugation space.
\end{Proposition}

\begin{proof}
Let $M=M_1\sharp M_2$
and let $N_i=M_i\setminus {\rm int}\Delta_i$.
By excision, one has ring isomorphisms
\begin{equation}\label{connsum-eq1}
H^*(M,N_1) \stackrel{\approx}{\to}
H^*(N_2,\partial N_2) \stackrel{\approx}{\leftarrow}
H^*(M_2,\Delta_2) \stackrel{\approx}{\to}
H^*(M_2,p_2) \, .
\end{equation}
The same isomorphisms hold for the $C$-equivariant cohomology and for the
cohomology of the fixed point sets. As $M_2$ is a \csp,
the pair $(M_2,p_2)$ is a \cpa\ by Remark~\ref{reduce}.
Therefore, $(M,N_1)$ is a \cpa.

\proref{EXPlac2} applied to $X=M_1$, $Y=N_1$ and $Z=\emptyset$
shows that $N_1$ is a \csp.
Applying then \proref{EXPlac}
to $X=M$, $Y=N_1$ and $Z=\emptyset$ proves
that $M$ is a \csp.  
\end{proof}

%%%%%%%%%%%%%%%%%%%%%%%%%%%%%%%%%%%%%%%%%%%%%%%%%%%%%%%%
\section{Conjugation complexes}\label{Sattcc}

\subsection{Attaching conjugation cells}\label{SSattcc}
Let $D^{2k}$ be the closed disk of radius $1$ in $\bbr^{2k}$, equipped
with an involution $\tau$ which is topologically conjugate to a
linear involution with exactly~$k$ eigenvalues equal to~$-1$.
As seen in Example~\ref{disk}, we call such a disk a {\it conjugation cell} of dimension~$2k$.
The fixed point set is then homeomorphic to a disk of dimension~$k$.
Observe that a product of two conjugation cells is a conjugation cell. 

Let $Y$ be a topological space with an involution $\tau$.
Let $\alpha\co S^{2k-1}\to Y$
be an equivariant map. Then the involutions on $Y$
and on $D^{2k}$ induce an involution on the quotient space
$$X= Y\cup_{\alpha} D^{2k}=
Y \dunion
D^{2k} \bigg/ \{u=\alpha(u)\mid x\in S^{2k-1}\}.
$$
We say that $X$ is obtained from $Y$ by attaching a {\it conjugation cell}
of dimension~$2k$. Note that the real locus $X^\tau$ is obtained
from $Y^\tau$ by adjunction of a $k$-cell.
Attaching a conjugation cell of dimension $0$ is making the disjoint union
with a point.

More generally, one can attach to $Y$ a set $\Lambda$ of $2k$-conjugation
cells, via an equivariant map
$\alpha\co \dunion_{\Lambda}S^{2k-1}_\lambda\to Y$.
The resulting space $X$ is equipped with an involution and its
real locus $X^\tau$ is obtained from $Y^\tau$ by adjunction of a collection
of $k$-cells labeled by the same set $\Lambda$.

The main result of this section is the following:

\begin{Proposition}\label{alcexten}
Let $(Y,Z)$ be a \cpa\ and let $X$ be obtained from $Y$ by attaching
a collection of conjugation cells of dimension $2k$.
Then $(X,Z)$ is a \cpa.
\end{Proposition}

\begin{proof} 
Without loss of generality, we may assume that
$Y$ and $X$ are path-connected.
We may also suppose that $Z$ and $Z^\tau$ are not empty.
Indeed, if $Z\neq\emptyset$ then $Z^\tau\neq\emptyset$
since $H^0(Y,Z)\approx H^0(Y^\tau,Z^\tau)$.
If $Z=\emptyset$, we replace
$Z$ by a point $pt\in Y^\tau$ ($Y^\tau$ is not empty if $Y$
is a \csp ) and use Remark~\ref{reduce}.

We shall now apply \proref{EXPlac}. The pair $(Y,Z)$ being
a \cpa\ by hypothesis, we must check that $(X,Y)$ is a \cpa.
By excision, $H^*(X,Y)=H^*(D,S)$, where
$D=\dunion_{\Lambda}D^{2k}_\lambda$ and
$S=\dunion_{\Lambda}S^{2k-1}_\lambda$.
One also has $H^*_C(X,Y)=H^*_C(D,S)$ and
$H^*(X^\tau,Y^\tau)=H^*(D^\tau,S^\tau)$, with
$D^\tau=\dunion_{\Lambda}D^{k}_\lambda$ and
$S^\tau=\dunion_{\Lambda}S^{k-1}_\lambda$.

Suppose first that $\Lambda=\{\lambda\}$ has one element,
so $D=D_\lambda$ and $S=S_\lambda$.
As seen in Example~\ref{disk}, we get here a \hfra\
$(\kappa_\lambda,\sigma_\lambda)$ such that,
if $a$ is the non-zero element of $H^{2k}(D,S)$, the equation
$r_\lambda\sigma_\lambda(a)=\kappa_\lambda(a)u^k$ holds.
For the general case, one has
$H^*(D,S)=\prod_{\lambda\in\Lambda}H^*(D_\lambda,S_\lambda)$,
$H^*_C(D,S)=\prod_{\lambda\in\Lambda}H^*_C(D_\lambda,S_\lambda)$,
etc, and
$\rho=\prod_{\lambda\in\Lambda}\rho_\lambda$,
$r=\prod_{\lambda\in\Lambda}r_\lambda$.
The homomorphisms $\sigma=\prod_{\lambda\in\Lambda}\sigma_\lambda$
and $\kappa=\prod_{\lambda\in\Lambda}\kappa_\lambda$ satisfy
$r\sigma(a)=\kappa(a)u^k$ for all $a\in H^{2k}(D,S)=H^*(D,S)$.
This shows that $(D,S)$ and then $(X,Y)$ is a \cpa.

We then know that $(X,Y)$ and $(Y,Z)$ are \cpa s. By \proref{EXPlac},
$(X,Z)$ is a \cpa. 
\end{proof}

\subsection{Conjugation complexes}\label{SCC}
Let $Y$ be a space with an involution $\tau$.
A space $X$ is a {\it \scc\ relative to } $Y$ if 
it is equipped with a filtration
$$
Y=X_{-1}\subset X_0\subset X_1\subset 
\cdots X={\scriptstyle \bigcup_{k=-1}^\infty} X_k \, 
$$
where $X_k$ is obtained from $X_{k-1}$ by the 
adjunction of a collection of conjugation cells
(indexed by a set $\Lambda_k(X)$).
The topology on $X$ is the direct limit topology of the $X_k$'s.
If $Y$ is empty, we say that $X$ is a {\it \scc}.
As in \cite{Gr}, the adjective ``spherical'' emphasizes 
that the collections of conjugation cells need not  
occur in increasing dimensions.

The involution $\tau$ on $Y$ extends naturally to an involution on $X$,
still called~$\tau$. The following result is a direct
consequence of \proref{alcexten} and \proref{Pcolim}.

\begin{Proposition}\label{SSCprop}
Let $X$ be a \scc\ relative to $Y$.
Then the pair $(X,Y)$ is a \cpa.  \qed
\end{Proposition}

\subsection{Remarks and Examples}\label{sccex}

\begin{ccote}\label{prodscc}\rm
Many topological properties of CW-complexes remain true for \scc es,
using minor adaptations of the standard techniques (see e.g.~\cite{LW}).
For instance, a \scc\ is paracompact, by the same proof as in~\cite[Theorem\,4.2]{LW}. 
Also, the product $X\times Y$ of two \csp s
admits a \scc-structure provided $X$ contains finitely
many conjugation cells, or both $X$ and $Y$ contain countably
many conjugation cells. For instance, one can order the elements
$(p,q)\in\bbn\times\bbn$ by the lexicographic ordering in $(p+q,p)$ and
construct a \csp\ $X\otimes Y$  by setting
$(X\times Y)_{(p,q)}=X_p\times X_q$. If $(p',q')$ is the successor of
$(p,q)$, then using that the product
of a conjugation cell is a conjugation cell, one shows that
$(X\times Y)_{(p',q')}$ is obtained from $(X\times Y)_{(p,q)}$
by adjunction of a collection of conjugation cells indexed by
$\Lambda_{(p',q')}(X\times Y)=\Lambda_{p'}(X)\times\Lambda_{q'}(Y)$.
There is then a $\tau$-equivariant continuous bijection
$\theta\co X\otimes Y\to X\times Y$. As in \cite[II.5,\,Theorem~5.2]{LW}, one shows that,
under the above hypotheses, $\theta$ is an homeomorphism.
\end{ccote}

\begin{ccote}\label{cpn-BT}\rm
The usual cell decomposition of $\bbc P^n$ ($n\leq\infty$)
makes the latter a \scc.
The product of finitely many copies of $\bbc P^\infty$
is also a \scc.  Here, we do not even need the preceeding
remark since we are just dealing with the product of countable CW-complexes.

Let $T$ be a torus (compact abelian group) of dimension $r$.
The involution $g\mapsto g^\mun$
induces an involution on the Milnor classifying space $BT$. The latter is
equivariantly homotopy equivalent to a product of $r$ copies of $\bbc P^\infty$
and therefore is a \csp. The isomorphism $\kappa_T$ of the \hfra\ for
$BT$ can be interpreted as follows.

Let $\hat T = {\rm Hom\,}(T,S^1)$ be the group of characters of $T$.
We have identifications
\begin{equation}\label{idenchi}
\hat T \approx [BT,\bbc P^\infty] \approx  H^2(BT;\bbz)  \ .
\end{equation}
Recall that $\hat T$ is a free abelian group of
rank the dimension of $T$. Hence
$H^{2}(BT)$ is isomorphic to
$\hat T\otimes\bbz_2$. 
For the $2$-torus subgroup $T_2$ of $T$, defined
to be the elements of $T$ of order $2$, 
one has in the same way
\begin{equation}{\rm Hom\,}(T_2,S^0)\approx  [BT_2,\bbr P^\infty]\approx
H^1(BT_2),\end{equation}
where we think of $S^0=\{\pm 1\}$ as the $2$-torus of $S^1$.
The homomorphism $\hat T\to {\rm Hom\,}(T_2,S^0)$, which sends $\chi\in\hat T$
to the restriction $\chi_{\scriptscriptstyle 2}$ of $\chi$ to $T_2$,
produces an isomorphism
$\kappa_{\scr T}\co H^{2}(BT)\to H^{1}(BT_2)$. Now the cohomology ring
$H^{2*}(BT)=S(H^{2}(BT))$ is the symmetric algebra over $H^{2}(BT)$,
and $H^{*}(BT_2)=S(H^{1}(BT_2))$. Therefore, the above isomorphism
$\kappa_{\scr T}$ extends to a ring isomorphism
$\kappa_{\scr T}\co H^{2*}(BT)\to H^{*}(BT_2)$
that is functorial in $T$.
Now, $BT_2=BT^\tau$,
and $\kappa_T$ is the isomorphism in the \hfra\ of $BT$.
This can be checked by choosing an isomorphism between $T$
and $(S^1)^r$, which induces a $C$-equivariant homotopy equivalence
between $BT$ and $(\bbc P^\infty)^r$ and a homotopy equivalence
between $BT^\tau$ and $(\bbr P^\infty)^r$.
\end{ccote}

\begin{ccote}\rm
Example~\ref{cpn-BT} generalizes to complex Grassmannians, with the
complex conjugation. The classical Schubert cells give the
\scc-structure. This generalizes to the coadjoint orbits of
compact semi-simple Lie groups with the Chevalley involution
(see Subsection~\ref{Cheva}), using the Bruhat-Schubert cells.
\end{ccote}

\begin{ccote}
Conjugation complexes with 3 conjugation cells.\ \rm
Let $X$ be a \scc\ with three conjugation cells, in dimension
$0$, $2k$ and $2l\geq 2k$. Then, $X$ is obtained by
attaching a conjugation cell $D^{2l}$ to the conjugation sphere
$\Sigma^{2k}$ (see Example~\ref{spheres}). 
The $C$-equivariant homotopy type of $X$
is determined by the class of the attaching map
$\alpha\in\pi^\tau_{2l-1}(\Sigma^{2k})$,
the equivariant homotopy group of $\Sigma^{2k}$ (the homotopy
classes of equivariant maps from $\Sigma^{2l-1}$ to $\Sigma^{2k}$).
We note $X=X_\alpha$.
Forgetting the $C$-equivariance and restricting to the fixed
point sets gives a homomorphism
$$\Phi_{l,k}\co \pi^\tau_{2l-1}(\Sigma^{2k})\to \pi_{2l-1}(S^{2k})\times
\pi_{l-1}(S^{k}).$$ In the case $k=1$ and $l=2$, this gives
$$\Phi\ccoeg \Phi_{2,1}\co \pi^\tau_{3}(\Sigma^2)\to \pi_3(S^2)\times\pi_1(S^1)
=\bbz\times\bbz.$$
\end{ccote}
Observe that the equivariant homotopy type of $X_\alpha$ and of
$X_\beta$ are distinct if $\Phi(\alpha)\neq\Phi(\beta)$. Indeed,
let $\Phi(\alpha)=(p,q)$. If  $a\in H^2(X;\bbz)$ and $b\in
H^4(X;\bbz)$ are the natural generators, then $a^2=pb$ 
(see, e.g.~\cite[\S\,9.5, Theorem\,3]{Sp}). 
Moreover $H^1(X^\tau;\bbz)=\bbz_q$. Note that since $X$ is a \csp, one has
$H^1(X^\tau)=\bbz_2$, which shows that $q$ must be even.

Now, it is easy to see that the Hopf map $h\co \Sigma^3\to\Sigma^2$
is $C$-equivariant; as $\Phi(h)=(1,2)$ is of infinite order,
this shows that there are infinitely many
$C$-equivariant homotopy types of \scc es with three conjugation
cells, in dimension $0$, $2$ and $4$.

%%%%%%%%%%%%%%%%%%%%%%%%%%%%%%%%%%%%%%%%%%%%%%%%%%%%%%%%%%%%%%%%%%%
\subsection{Equivariant fiber bundles over \scc es}\label{Sfibscc}

Let $G$ be a topological group together with an involution $\sigma$
which is an automorphism of $G$. Let $(B,\tau)$ be a space with involution.
By a {\it $(\sigma,G)$-principal bundle} we mean a (locally trivial)
$G$-principal bundle $p\co E\to B$ together with an involution
$\tilde\tau$ on $E$ satisfying $p\pcirc\tilde\tau=\tau\pcirc p$ and
$\tilde\tau(z\cdot g)=\tilde\tau(z)\cdot\sigma(g)$ for all $z\in E$ and
$g\in G$. Following the terminology of \cite[p.\,56]{tD}, 
a $(\sigma,G)$-principal bundle is a $(C,\check\sigma,G)$-bundle, where
$\check\sigma\co C\to G$ is the homomorphism sending the generator of $C$
to $\sigma$.

Let $F$ be a space together with an involution $\tau$ and a left
$G$-action. We say that the involution $\tau$ and
the $G$-action are {\it compatible} if $\tau(gy)=\sigma(g)\,\tau(y)$.
This means that the $G$-action extends to an action of
the semi-direct product $G^\times=G\rtimes C$.

Let $p\co E\to B$ be a $(\sigma,G)$-principal bundle.
Let $(F,\tau)$ be a space with involution together with a compatible $G$-action.
The space $E\times_G F$ inherits an involution (also called $\tau$)
and the associated bundle $E\times_G F\to B$, with fiber $F$,
is a $\tau$-equivariant locally trivial bundle.

\begin{Proposition}\label{Pfibscc}
Suppose that $G$ is a compact Lie group, that $F$
is a \csp\ and that $B$ is a \scc. Then $E\times_G F$
is a \csp.
\end{Proposition}

\begin{proof}
Suppose first that $B=D$ is a conjugation cell, with boundary $S$.
Then $E$ is compact and, by \cite[Ch.\,1, Proposition~8.10]{tD}, $p$
is a {\it locally trivial $(C,\tilde\sigma,G)$-bundle}.
This means that there exists an open covering $\calu$ of $B$ by
$C$-invariant sets such that for each $U\in\calu$ the bundle
$p^\mun(U)\to U$ is induced by a $(\sigma,G)$-principal bundle
over a $C$-orbit (namely one point or two points). Since the quotient space
$C\backslash D$ is compact, the coverings $\calu$ admits a
partition of unity by $C$-invariant maps. Together these imply that
the $(\sigma,G)$-bundle $p$ is induced from a universal 
$(C,\tilde\sigma,G)$-bundle by a
$C$-equivariant map from $D$ to some classifying space  
and $C$-homotopic maps induce isomorphic $(C,\tilde\sigma,G)$-bundles
\cite[Ch.\,1, Theorem\,8.12 and 8.15]{tD}. The cell $D$ is
$C$-contractible, which implies that $E=D\times G$
and $E\times_G F=D\times F$, with the product involution.
By \proref{Pprod}, the pair $(E,E_{|S})$ is a \cpa.

This enables us to prove \proref{Pfibscc} by induction on
the $n$-stage $B_n$ of the construction of $B$ as a \scc.
Let $Z_n=p^\mun(B_n)\times_G F$. As $B_0$ is discrete, $Z_0$ is
the disjoint union of copies of $F$ and is then a \csp.
Suppose by induction that $Z_{n-1}$ is a \csp.
The above argument shows that $(Z_n,Z_{n-1})$ is a \cpa.
Using \proref{EXPlac}, one deduces that $Z_n$ is a \csp.
Therefore, $Z_n$ is a \csp\ for all $n\in\bbn$.
By \proref{Pcolim}, this implies that $E\times_G F=\lim_\to Z_n$ is a \csp.
\end{proof}

\begin{Remark}\label{Pfibscc-rem}
An analogous argument also gives a relative version of
\proref{Pfibscc} for pairs of bundles over $X$,
with a \cpa\ of fibers $(F,F_0)$. The same remains true for
a bundle over a relative \scc.
\end{Remark}

%%%%%%%%%%%%%%%%%%%%%%%%%%%%%%%%%%%%%%%%%%%%%%%%%%%%%%%%%%%%%%%%%%%
\section{Conjugate-equivariant complex bundles}\label{tabdles}

%%%%%%%%%%%%%%%%%%%%%%%%%%%%%%%%%%%%%%%%%%
\subsection{Definitions}
Let $(X,\tau)$ be a space with an involution. A {\it
$\tau$-conjugate-equivariant bundle} (or, briefly, a {\it \taub})
over $X$ is a complex vector bundle $\eta$, with total space
$E=\bbe(\eta)$ and bundle projection $p\co E\to X$, together with an
involution $\hat\tau\co E\to E$ such that $p\pcirc\hat\tau=\tau\pcirc
p$ and $\hat\tau$ is conjugate-linear on each fiber:
$\hat\tau(\lambda\,x)=\bar \lambda\hat\tau(x)$ for all
$\lambda\in\bbc$ and $x\in E$. Atiyah was the first to study \taub
s \cite{At}.  He called them ``real bundles'' and used them to
define $KR$-theory.

Let $P\to X$ be a $(\sigma,U(r))$-principal bundle in the sense of
Subsection~\ref{Sfibscc}, with $\sigma\co U(r)\to U(r)$ being the complex
conjugation. Then, the associated bundle $P\times_{U(r)}\bbc^r$,
with $\bbc^r$ equipped with the complex conjugation, is a \taub\
and any \taub\ is of this form. It follows that if  
$p\co E\to X$ be a \taub\ $\eta$ of rank $r$ and if
$E^{\hat\tau}$ is the fixed point set of $\hat\tau$, then
$p\co E^{\hat\tau}\to X^\tau$ is a real vector bundle
$\eta^\tau$ of rank $r$ over $X^\tau$. 

Examples of \taub\ include the canonical complex vector bundle over
$BU(r)$ or over the complex Grassmannians. Note that a bundle induced
from a \taub\ by a $C$-equivariant map is a \taub. 

\begin{Proposition}\label{taubind}
Let $\eta$ be a \taub\ of rank $r$ over a space with involution $(X,\tau)$.
If $X$ is paracompact, then $\eta$ is induced from the universal
bundle by a $C$-equivariant map from $X$ into $BU(r)$. Moreover,
two $C$-equivariant map which are $C$-homotopic induce isomorphic \taub s.
\end{Proposition}

\begin{proof} 
It is equivalent to prove the corresponding statement of \proref{taubind} 
for $(\sigma,U(r))$-bundles.
Let $p\co P\to X$ be a $(\sigma,U(r))$-bundle. As $X$ is paracompact and $U(r)$
is compact, the total space $P$ is paracompact. Therefore, by 
\cite[Ch.\,1, Proposition~8.10]{tD}, $p$
is a locally trivial $(\sigma,U(r))$-bundle,
meaning that there exists an open covering $\calv$ of $X$ by
$C$-invariant sets such that for each $V\in\calv$ the bundle
$p^\mun(V)\to V$ is induced by a $(\sigma,G)$-principal bundle 
$q\co q_\calo\to\calo$ over a $C$-orbit $\calo$. 
When $\calo$ consists of one point $a$, one can identify $Q_\calo$ with $U(r)$ 
such $\tilde\tau(\gamma)=\bar\gamma$. 
For a free orbit $\calo=\{a,b\}$,
one can identify $Q_\calo$ with $\calo\times U(r)$ such that
$\tilde\tau(a,\gamma)=(b,\bar\gamma)$ and $\tilde\tau(b,\gamma)=(a,\bar\gamma)$.
Using these, one gets a family of $U(r)$-equivariant maps 
$\{\varphi_{\scr V}\co p^\mun(V)\to U(r)\mid V\in\calv\}$ such that
\begin{equation}\label{taubind-eq1}
\varphi_{\scr V}\pcirc\tau(z)=\overline{\varphi_{\scr V}(z)} \, , 
\end{equation}
for all $V\in\calv$. The quotient space
$C\backslash X$ is also paracompact. Therefore, the coverings $\calv$ admits a
locally finite partition of the unity $\mu_{\scr V}$, $V\in\calv$, by $C$-invariant maps.
Using $\{\varphi_{\scr V},\mu_{\scr V}\mid V\in\calv\}$, we can perform the classical
Milnor construction of a map $f\co X\to BU(r)$ inducing $p$. Because of
Equation~\eqref{taubind-eq1}, $f$ is $C$-equivariant. 
The last statement of \proref{taubind} is a direct consequence of  
\cite[Ch.\,1, Theorem\,8.12 and 8.15]{tD}. 
\end{proof}

\begin{Corollary}\label{taubccell}
Let $\eta$ be a \taub\ over a conjugation cell. Then, the total space
of disk bundle $\bbd(\eta)$ is a conjugation cell.
\end{Corollary}

\begin{proof} 
As a conjugation cell is $C$-contractible,  \proref{taubind}
implies that $\eta$ is a product bundle.  We then use that the product 
of two conjugation cells is a conjugation cell. 
\end{proof}

\begin{Remark} Pursuing in the way of \proref{taubind}, one can 
prove that the set of isomorphism classes of \taub s of rank $r$ over a
paracompact space $X$ is in bijection with the set of $C$-equivariant
homotopy classes of $C$-equivariant maps from $X$ to $BU(r)$. 
\end{Remark}

%%%%%%%%%%%%%%%%%%%%%%%%%%%%%%%%%%%%%%%%%%%%%%%%%%%%%%%%
\subsection{Thom spaces}

\begin{Proposition}\label{thomcbdle}
Let $\eta$ be a \taub\
over a \csp\ $X$. Then the total space
$\bbd(\eta)$ of the disk bundle of $\eta$ and the total space
$\bbs(\eta)$ of the sphere bundle of $\eta$ form a \cpa\
$(\bbd(\eta), \bbs(\eta))$.
\end{Proposition}

\proof
Let $\bbe(\eta)\to X$ be the bundle projection and let
$r$ be the rank of $\eta$. Performing the
Borel construction $\bbe(\eta)_C\to X_C$ gives a complex
bundle $\eta_{\scr C}$ of rank $r$  over $X_C$ and $\eta$ is induced from
$\eta_{\scr C}$ by the map $X\to X_C$.
The following diagrams, in which
the letters $\calt$ denote the Thom isomorphisms,
show how to define $\overline{\sigma}$ and $\overline{\kappa}$.
$$
\xymatrix@C-3pt@M+2pt@R-4pt{%
H^{2*-2r}_C(X) \ar[d]_(0.50){\rho}
\ar[r]^(0.42){\calt_C}_(0.42){\approx}  &
H^{2*}_C(\bbd(\eta),\bbs(\eta)) \ar[d]_(0.50){\rho}  \\
%%%%%% ROW 2
H^{2*-2r}(X) \ar@/_/[u]_{\sigma}
\ar[r]^(0.42){\calt}_(0.42){\approx}  &
H^{2*}(\bbd(\eta),\bbs(\eta))
\ar@/_/[u]_{\bar\sigma\co =\calt_{\scr C}\pcirc\sigma\pcirc\calt^\mun}
}
$$
$$
\xymatrix@C-3pt@M+2pt@R-4pt{%
H^{2*-2r}(X) \ar[d]_(0.50){\kappa}
\ar[r]^(0.42){\calt}_(0.42){\approx}  &
H^{2*}(\bbd(\eta),\bbs(\eta)) \ar[d]^(0.50){\bar\kappa\co =\calt^\tau\pcirc\kappa\pcirc\calt^\mun}  \\
%%%%%% ROW 2
H^{*-r}(X^\tau)
\ar[r]^(0.42){\calt^\tau}_(0.42){\approx}  &
H^{*}(\bbd(\eta^\tau),\bbs(\eta^\tau))
}
$$

Consider also the following commutative diagram,
where the vertical arrows are restriction to a fiber.
$$
\xymatrix{
H^{2r}(\bbd(\eta),\bbs(\eta)) \ar[r]^{\bar\sigma} \ar[d]_{j} &
H^{2r}_C(\bbd(\eta),\bbs(\eta)) \ar[r]^{\overline{r}} \ar[d]_{\bar j} &
H^{2r}_C(\bbd(\eta)^\tau,\bbs(\eta)^\tau) \ar[d]^{j^\tau}\\
H^{2r}(D^{2r},S^{2r-1}) \ar[r]^{\sigma_{\scr D}} &
H_C^{2r}(D^{2r},S^{2r-1}) \ar[r]^{r_{\scr D}} & H_C^{2r}(D^r,S^{r-1})
}$$
It remains to prove the conjugation equation.
Let $\tho(\eta)\in H^{2r}(\bbd(\eta),\bbs(\eta))$ be the Thom class
of $\eta$.
By definition of $\bar\sigma$,
one has $\bar\sigma(\tho(\eta))=\tho(\eta_{\scr C})$.
Observe that $D^{2r}$ is a conjugation cell
and $\bar j (\tho(\eta_{\scr C}))=\sigma_{\scr D}([D^{2r},S^{2r-1}])$.
Therefore
\begin{equation}
r\pcirc\sigma_{\scr D}\pcirc j (\tho(\eta)) =
\kappa_{\scr D^{2r}}([D^{2r},S^{2r-1}])\, u^r
= [D^{r},S^{r-1}]\, u^r
\ .
\end{equation}
But $r\pcirc\sigma_{\scr D}\pcirc j=j^\tau\pcirc \bar r\pcirc\bar\sigma$
and the preimage under $j^\tau$ of $[(D^{r},S^{r-1})]$
is $\tho (\eta^\tau)$. By \lemref{hfpri}, the kernel of
$j^\tau H^{2r}_C(\bbd(\eta)^\tau,\bbs(\eta)^\tau)\to
H_C^{2r}(D^r,S^{r-1})$ is of type $\lt{r}$. Therefore, one has
\begin{equation}\label{thomcbdle-eq2}
\bar r\pcirc\bar\sigma (\tho(\eta)) = \bar r(\tho(\eta_{\scr C})) =
\tho(\eta^\tau)\, u^r + \lt{r}
\ .
\end{equation}
Using Equation~\eqref{thomcbdle-eq2}, one has, for $x\in H^{2k+2r}(\bbd(\eta),\bbs(\eta))$:
$$
\begin{array}{rcl}
\bar r\pcirc\bar\sigma(x) &=&
\bar r\pcirc\calt_{\scr C}\pcirc\sigma\pcirc\calt^\mun(x)  =
\bar r\big(\tho(\eta_{\scr C})\cdot \sigma\pcirc\calt^\mun(x)\big)
\\[2mm] &=&
\bar r \big(\tho(\eta_{\scr C})\cdot r\pcirc\sigma\pcirc\calt^\mun(x)\big) \\[2mm] &=&
(\tho(\eta^\tau)\, u^r + \lt{r})\,(\kappa(\calt^\mun(x))u^k+\lt{k})
\\[2mm] &=&
\tho(\eta^\tau)\,\kappa(\calt^\mun(x))\, u^{k+r} +\lt{k+r} =
\bar\kappa(x)\, u^{k+r} +\lt{k+r} \, . 
\end{array}\eqno{\smash{\lower 27pt\hbox{\qed}}}
$$

\begin{Remark}\rm
The pair $(\bbd(\eta),\bbs(\eta))$ is cohomologically equivalent to the
pair $(\bbd(\eta)/\bbs(\eta),pt)$ and $\bbd(\eta)/\bbs(\eta)$
is the Thom space of $\eta$. Using Remark~\ref{reduce},
\proref{thomcbdle} says that
if $\eta$ is a \taub\
over a conjugation space, then the Thom space of $\eta$ is a \csp.
\end{Remark}

\begin{Remark}\label{eulercl}\rm
By the definition of $\bar\kappa\co H^{2*}(\bbd(\eta),\bbs(\eta))\to
H^*(\bbd(\eta^\tau),\bbs(\eta^\tau))$, one has $\bar\kappa(\tho(\eta))=
\tho(\eta^\tau))$. The inclusion $(\bbd(\eta),\emptyset)\subset(\bbd(\eta),\bbs(\eta))$
is a $C$-equivariant map between conjugation pairs and $\bbd(\eta)$ is
$C$-homotopy equivalent to $X$. The induced homomorphisms on cohomology
$i\co H^{2r}(\bbd(\eta),\bbs(\eta))\to H^{2r}(X)$ and
$i^\tau\co H^{r}(\bbd(\eta^\tau),\bbs(\eta^\tau))\to H^{r}(X^\tau)$
send the Thom classes $\tho(\eta)$ and $\tho(\eta^\tau)$ to the Euler classes
$e(\eta)$ and $e(\eta^\tau)$. By naturality of the \hfra s, we deduce that,
for any \ceb\ $\eta$ over a \csp\ $X$, one has $\kappa(e(\eta))=e(\eta^\tau)$.
This will be generalized in \proref{chernSW}.
\end{Remark}

We finish this subsection with the analogue
of \proref{thomcbdle} for \scc es.

\begin{Proposition}\label{ThomSCC}
Let $\eta$ be a $\tau$-bundle
over a \scc\ $X$. Then,
$\bbd(\eta)$  is a  \scc\ relative to $\bbs(\eta)$.
\end{Proposition}

\begin{proof}
Let $X$ be obtained from $Y$ by attaching
a collection of conjugation cells of dimension $2k$,
indexed by a set $\Lambda$.  Let
$D=\dunion_{\Lambda}D^{2k}_\lambda$ and
$S=\dunion_{\Lambda}S^{2k-1}_\lambda$ ($\lambda\in\Lambda$).
Let $\pi=\pi_{\scr D}\dunion \pi_{\scr Y} \co  D\dunion Y \to X$
be the natural projection. Then
$\bbd(\eta)$ is obtained from $\bbd(\pi_{\scr Y}^*\eta)\cup\bbs(\eta)$
by attaching $\bbd(\pi_{\scr D}^*\eta)$.  By Corollary~\ref{taubccell},
$\bbd(\pi_{\scr D_\lambda}^*\eta)$ is a conjugation cell of dimension
$2k+2r$, where $r$ is the complex rank of $\eta$.
Therefore, $\bbd(\eta)$ is obtained from
$\bbd(\pi_{\scr Y}^*\eta)\cup\bbs(\eta)$
by attaching a collection of conjugation cells of dimension $2k+2r$.
This proves \proref{ThomSCC}. 
\end{proof}

%%%%%%%%%%%%%%%%%%%%%%%%%%%%%%%%%%%%%%%%%%%%%%%%%%%%%%%%
\subsection{Characteristic classes}\label{split}

If $\eta$ be a \taub\ over a space with involution
$X$, we denote by $c(\eta)\in H^{2*}(X)$ the ($\mod\ 2$) total Chern class
of $\eta$ and by $w(\eta^\tau)\in H^{*}(X^\tau)$ the
total Stiefel-Whitney class of $\eta^\tau$. The aim of this section is to
prove the following:

\begin{Proposition}\label{chernSW}
Let $\eta$ be a $\tau$-bundle
over a \scc\ $X$. Then
$\kappa(c(\eta))=w(\eta^\tau)$.
\end{Proposition}

\begin{proof}
Let $q\co  \bbp(\eta)\to X$ be the projective bundle associated to $\eta$,
with fiber $\bbc P^{r-1}$. The conjugate-linear involution $\hat\tau$
on $\bbe(\eta)$ descends to an involution $\tilde\tau$ on $\bbp(\eta)$
for which the projection $q$ is equivariant.
One has $\bbp(\eta)^{\tilde\tau}=\bbp(\eta^\tau)$,
the projective bundle associated to $\eta^\tau$, with fiber $\bbr P^{r-1}$.
We also call $q\co \bbp(\eta^\tau)\to X^\tau$
the restriction of $q$ to $\bbp(\eta^\tau)$.

As $q$ is equivariant, the induced complex vector bundle $q^*\eta$
is a $\tilde\tau$-bundle with
$\bbe(q^*\eta)^\tau=\bbe(q^*\eta^\tau)$.
Recall that $q^*\eta$ admits
a canonical line subbundle $\lambda_\eta$:  a point of
$\bbe(\lambda_\eta)$ is a couple $(L,v)\in \bbp(\eta)\times\bbe(\eta)$
with $v\in L$. The same formula holds for $\eta^\tau$, giving a
real line subbundle $\lambda_{\eta^\tau}$ of $q^*\eta^\tau$.
Moreover, $\hat\tau(v)\in\tau(L)$ and thus
$\lambda_\eta$ is a $\tilde\tau$-conjugate-equivariant line bundle
over $\bbp(\eta)$. Again, $\bbe(\lambda_\eta)^\tau=\bbe(\lambda_{\eta^\tau})$.
The quotient bundle $\eta_1$ of $\eta$ by $\lambda_\eta$
is also a $\tilde\tau$-bundle over $\bbp(\eta)$
and $q^*\eta$ is isomorphic to the equivariant Whitney sum of
$\lambda_\eta$ and $\eta_1$.

By \proref{Pfibscc}, $\bbp(\eta)$ is a \csp.
Denote by $(\tilde\kappa,\tilde\sigma)$
its \hfra. By \remref{eulercl}, one has $\tilde\kappa(c_1(\lambda_\eta))=w_1(\lambda_{\eta^\tau})$.
As $\tilde\kappa$ is a ring isomorphism,
one has $\tilde\kappa(c_1(\lambda_\eta)^k)=w_1(\lambda_{\eta^\tau})^k$ for
each integer $k$.

By \cite[Chapter 16,2.6]{Hu}, we have
in $H^{2*}(\bbp(\eta))$ the equation
\begin{equation}\label{chern-eq}
c_1(\lambda_\eta)^r=
\sum_{i=1}^{r} q^*(c_i(\eta))\, c_1(\lambda_\eta)^{r-i} .
\end{equation}
and, in $H^{*}(\bbp(\eta^\tau))$,
\begin{equation}\label{chern-eqtau}
w_1(\lambda_{\eta^\tau})^r=
\sum_{i=1}^{r} q^*(w_i(\eta^\tau))\, w_1(\lambda_{\eta^\tau})^{r-i} \, .
\end{equation}
As $\tilde\kappa(c_1(\lambda_\eta))=w_1(\lambda_{\eta^\tau})$
and $\tilde\kappa\pcirc q^* = q^*\pcirc\kappa$, applying
$\tilde\kappa$ to Equation~\eqref{chern-eq} and using
Equation~\eqref{chern-eqtau} gives
\begin{equation}\label{chern-eqfin}
\sum_{i=1}^{r} q^*(\kappa(c_i(\eta)))\, w_1(\lambda_{\eta^\tau})^{r-i} =
\sum_{i=1}^{r} q^*(w_i(\eta^\tau))\, w_1(\lambda_{\eta^\tau})^{r-i} \, .
\end{equation}
By the Leray-Hirsch theorem,
$H^*(\bbp(\eta^\tau))$ is a free $H^*(X^\tau)$-module with basis
$w_1(\lambda_\eta)^k$ for $k=1,\dots,r-1$,
and $q^*$ is injective.
Therefore, Equation~\eqref{chern-eqfin} implies
\proref{chernSW}. 
\end{proof}

\begin{Remark}\rm By \proref{taubind}, it would be enough to
prove \proref{chernSW} for the canonical bundle over the Grassmannian.
This can be done via the Schubert calculus
(see \cite[Problem 4-D, p. 171, and \S\,6]{MS}). Such an argument proves
\proref{chernSW} for $X$ a paracompact \csp.
\end{Remark}

\section{Compatible torus actions}\label{Sctac}

Let $X$ be a space together with an involution $\tau$. Suppose
that a torus $T$ acts continuously on $X$. We say that the
involution $\tau$ is {\it compatible} with this torus action if
$\tau(g\cdot x)=g^\mun\cdot\tau(x)$ for all $g\in T$ and $x\in X$.
It follows that $\tau$ induces an involution on on the fixed point
set $X^T$. Moreover, the $2$-torus subgroup $T_2$ of $T$, defined
to be the elements of $T$ of order $2$, acts on $X^\tau$.
The involution and the $T$-action extend to an action of the
semi-direct product $T^\times=T\rtimes C$, where $C$ acts on $T$ by 
$\tau\cdot g=g^\mun$.

When a group $H$ acts on $X$, we denote by $X_H$ the Borel construction
of $X$. Observe that if $T^\times$ as above acts on $X$, then the diagonal
action of $C$ on $ET\times X$ descends to an action of $C$ on $X_T$.

\begin{Lemma}\label{borbor}
Let $X$ be a space together with a continuous action of $T^\times$. Then
$X_{T^\times}$ has the homotopy type of $(X_T)_C$.
\end{Lemma}

\begin{proof} 
$X_T$ has the $C$-equivariant homotopy type of the quotient
$T\backslash ET^\times\times X$, where $T$ acts on $ET^\times\times X$ by
$g\cdot(w,x)=(wg^\mun,gx)$. The formula
$\tau\cdot(w,x)=(w\tau,\tau(x))$ then induces a $C$-action on
$X_T$ which is free. Therefore, $(X_T)_{T^\times}=C\backslash X_T=X_{T^\times}$.
\end{proof}

The particular case of $X=pt$ in Lemma~\ref{borbor} gives the following:

\begin{Corollary}\label{BGBT}
$BT^\times\simeq BT_C$. \qed
\end{Corollary}

\begin{Lemma}\label{CBTfix}
$(X_T)^\tau=(X^\tau)_{T_2}$.
\end{Lemma}

\begin{proof} 
Let $H$ be a group acting continuously on a space $Y$.
Recall that elements of the infinite joint $EH$ are represented by
sequences $(t_ih_i)$ ($i\in\bbn$) with $h_i\in H$ and  $t_i\in
[0,1]$, almost all vanishing, with $\sum t_i=1$. Under the right
diagonal action of $H$ on $EH$, each $(t_ih_i)$ is equivalent to a
unique element $(t_i\tilde h_i)$ for which $\tilde h_j=I$, the
unit element of $H$, where $j$ is the minimal integer $k$ for
which $t_k\neq 0$. Therefore, each class in $BH=EH/H$ has a unique
such representative which we call {\it minimal}. In the same way,
each class in $Y_H$ has a unique {\it minimal representative}
$(w,y)\in EH\times Y$ for which $w$ is minimal.

One easily check that there is a commutative diagram:
\begin{equation}\label{CBTfix-eq1}
\xymatrix@C-3pt@M+2pt@R-2pt{%
(X^\tau)_{T_2} \ar[drr]_(0.55){\beta}
\ar@{>->}[r]  &
X_{T_2} \ar@{>>}[r] &
X_T   \\
%%%%%% ROW 2
&&  (X_T)^\tau \ar@{>->}[u]
}
\end{equation}
Working with minimal representatives in $(X^\tau)_{T_2}$, we see
that the natural map $(X^\tau)_{T_2}\to X_T$ is injective.
Hence, $\beta$ is injective. Let $(w,x)\in ET\times X$
with $w=(t_iz_i)$ minimal. Then, $\tau(w,x)$ is also a minimal
representative. If $\tau(w,x)=(w,x)$ in $X_T$, this
implies that $\tau(x)=x$ and
$z_i^\mun=z_i$, that is $z_i\in T_2$ (when $t_i\neq 0$).
This proves that $\beta$ is surjective.  
\end{proof}

\begin{Example}\label{XS1g2}\rm
Let $X=S^1\subset\bbc$ with the complex conjugation as involution,
and $T=S^1$ acting on $X$ by $g\cdot z = g^2z$. Then, $X^\tau=S^0$
on which $T_2$ acts trivially, so $(X^\tau)_{T_2} = BT_2\times S^0$.
On the other hand, $X$ is a $T$-orbit so $X_T=ET/T_2$.
The space $ET/T_2$ has the homotopy type of $BT_2$ but
$(ET/T_2)^\tau$ has two connected components, both
homeomorphic to $BT_2$. One is the image
of $(ET)^\tau=ET_2$ and is equal to $\beta(BT_2\times\{1\})$.
The other is the image of $\{(t_jh_j)\mid h_j=\pm i\}$ and is
equal to $\beta(BT_2\times\{-1\})$.
\end{Example}

The main result of this section is the following:

\begin{Theorem}\label{strcomp}
Let $(X,Y)$ be a \cpa\ together with a compatible action of
a torus $T$. Then, the involution induced on $(X_T,Y_T)$ is a conjugation.
\end{Theorem}

\begin{proof} 
Assume first that $Y=\emptyset$.
The universal bundle $p\co ET\to BT$ is a $(T,\sigma)$-principal bundle
in the sense of Subsection~\ref{Sfibscc}, with $\sigma(g)=g^\mun$,
and $X_T\to BT$ is the
associated bundle with fiber $X$.
As $BT$ is a \csp\ (see Remark~\ref{cpn-BT}
in Subsection~\ref{sccex}), the space $X_T$ is a \csp\ by
\proref{Pfibscc}. When $Y$ is not empty, we use Remark~\ref{Pfibscc-rem}.
\end{proof}

Using Lemma~\ref{CBTfix}, one gets the following corollary of
Theorem~\ref{strcomp}.

\begin{Corollary}\label{srtcomp-cor}
Let $X$ be a space together with an involution and a compatible
$T$-action. Then, there is a ring isomorphism
$$\bar\kappa\co  H^{2*}_T(X)\hfl{\approx}{} H^{*}_{T_2}(X^\tau). \eqno{\qed}$$
\end{Corollary}

We end this section with a result that will be used in Section~\ref{SHam}.
Let $\eta$ be a $T$-equivariant $\tau$-bundle over a space with
involution $X$. Precisely, $\eta$ is a $\tau$-bundle over $X$ and
there is a $\hat\tau$-compatible $T$-action on $\bbe(\eta)$,
over the identity of $X$, which is $\bbc$-linear on each fiber.
Let $r$ be the complex rank of $\eta$.
The $T$-Borel construction on $\bbe(\eta)\to X$
produces a complex vector bundle  $\eta_{\scr T}$
of rank $r$ over $X_T$. One checks that the involution induced on
$\bbe(\eta_{\scr T})=\bbe(\eta)_{\scr T}$ makes
$\eta_{\scr T}$ a $\tau$-bundle (the letter $\tau$ also
denotes here the involution induced on $X_T=BT\times X$).
For a $T^\times$-invariant Riemannian metric on $\eta$, the spaces
$\bbd(\eta)$ and $\bbs(\eta)$ are $T^\times$-invariant.

\begin{Proposition}\label{thomcbdleTeq}
Let $\eta$ be a $T$-equivariant $\tau$-bundle over
a \csp\ $X$. Then the pair $(\bbd(\eta)_T,\bbs(\eta)_T)$
is a \csp.
\end{Proposition}

\begin{proof}
As the Riemannian metric is $T^\times$-invariant,
one has $\bbd(\eta)_{\scr T}=\bbd(\eta_{\scr T})$
and $\bbs(\eta)_{\scr T}=\bbs(\eta_{\scr T})$. By
Theorem~\ref{strcomp}, the base space $BT\times X$
of $\eta_{\scr T}$ is a \csp.
\proref{thomcbdleTeq} then follows from
\proref{thomcbdle}.  
\end{proof}

%%%%%%%%%%%%%%%%%%%%%%%%%%%%%%%%%%%%%%%%%%%%%%%%%%%%%%%%%%%%%%%%%%%
\section{Hamiltonian manifolds with anti-symplectic involutions}\label{SHam}
\subsection{Preliminaries}\label{Hampreli}

Let $M$ be a compact symplectic manifold equipped with a
Hamiltonian action of a torus $T$. Let $\tau$ be a smooth
anti-symplectic involution on $M$ compatible with the action of
$T$ (see Section~\ref{Sctac}). Thus, the semi-direct group 
$T^\times\ccoeg T\rtimes C$ acts on $M$. Moreover, if it is
non-empty, $M^\tau$ is a Lagrangian submanifold, called the {\em
real locus} of $M$. For general work on such involutions
together with a Hamiltonian group action, see \cite{Du} and
\cite{OS}.

We know that the symplectic manifold $(M,\omega)$ admits an almost
Kaehler structure calibrated by $\omega$. That is, there is an
almost complex structure $J\in{\rm End\,}TM$ together with a
Hermitian metric $h$ whose imaginary part is $\omega$ (see
\cite[\S\,1.5]{Au}; $J$ and $h$ determine each other). These
structures form a convex set and by averaging, we can find an
almost complex structure whose Hermitian metric $\tilde h$ is
$T$-invariant. Now, the Hermitian metric
\begin{equation}h(v,w)\ccoeg \frac{1}{2}\Big(\tilde h(v,w) +
\overline{\tilde h((T\tau(v),T\tau(w))}\Big)\end{equation}
is still $T$-invariant and satisfies
$h(T\tau(v),T\tau(w))=\overline{h(v,w)}$.
We suppose that the symplectic manifold
$(M,\omega)$ is equipped with such an  almost
Kaehler structure $(J,h)$ calibrated by $\omega$,
which we call a {\it $T^\times$-invariant almost Kaehler structure}.

Let $\Phi\co M\to\algt^*$ be a moment map for the Hamiltonian torus
action, where $\algt$ denotes the Lie algebra of $T$
and $\algt^*$ denotes its vector space dual. 
Evaluating $\Phi$ on a generic element $\xi$ of $\algt$ yields a real
Morse-Bott function $\Phi^\xi(x)=\Phi(x)(\xi)$ whose critical
point set is $M^T$.  Suppose $F$ is a connected component of
$M^T$. By \cite[\S\,III.1.2]{Au}, $F$ is an almost Kaehler (in
particular symplectic) submanifold of $M$. If
$F^\tau\neq\emptyset$, then $F$ is preserved by $\tau$:
$\tau(F)=F$.

Let $\nu(F)$ be the normal bundle to $F$, seen as the orthogonal
complement of $TF$. The bundle $\nu(F)$ is then a complex vector
bundle. By $T^\times$-invariance of the Hermitian metric, $\nu(F)$ admits
a $\bbc$-linear $T$-action and $\tau\co F\to F$ is covered by an
$\bbr$-linear involution $\hat\tau$ of the total space
$\bbe(\nu(F))$ which is compatible with the $T$-action. Moreover,
$\nu(F)$ inherits a Hermitian metric $h$ whose imaginary part is
the symplectic form $\omega$. Let $x\in F$. For
$v\in\bbe_x(\nu(F))$, $w\in\bbe_{\tau(x)}(\nu(F))$ and
$\lambda\in\bbc$, one has
\begin{equation}\label{etasq}
\begin{array}{rcl}
h(\hat\tau(\lambda v),w) & = &\overline{h(\lambda v,\hat\tau(w))}
=\bar\lambda\,\overline{h(v,\hat\tau(w))}\\
& = & \bar\lambda\,h(\hat\tau(v),w) =h(\bar\lambda\hat\tau(v),w).
\end{array}
\end{equation}
This shows that $\nu(F)$ is
a $\tau$-bundle.

Let us decompose $\nu(F)$ into a Whitney sum of $\chi$-weight
bundles $\nu^\chi(F)$ for $\chi\in\hat T$, the group of smooth homomorphisms
from $T$ to $S^1$. Recall that the latter is free abelian of rank the
dimension of $T$. We call $\nu^\chi(F)$
an {\it isotropy weight bundle}. Since the $T$-action on
$\nu(M^T)$ is compatible with $\hat\tau$, the isotropy weight
bundles are preserved by $\hat\tau$ and are thus $\tau$-bundles.
Consequently, the {\it negative normal bundle} $\nu^-(F)$,
which is the Whitney sum of those $\nu^\chi(F)$ for which
$\Phi^\xi(\chi)<0$, is a $\tau$-bundle.

Of course $M^T\subset M^{T_2}$. The case where this inclusion
is an equality will be of interest.
\begin{Lemma}\label{nuchi}
The following conditions are equivalent:
\begin{enumerate}
\renewcommand{\labelenumi}{\rm(\roman{enumi})}
\item\label{nuchi2} $M^{T}=M^{T_2}$.
\item\label{nuchi1} $M^\tau\cap M^{T}=(M^\tau)^{T_2}$.
\item\label{nuchi3} for each $x\in M^T$,
there is no non-zero weight $\chi\in\hat T$ of the
isotropy representation of $T$ at $x$ such that $\chi\in 2\cdot\hat T$.
\end{enumerate}
\end{Lemma}

\begin{proof}
If (ii) is true, then
\begin{equation}
M^\tau\cap M^{T}\subset (M^\tau)^{T_2} = M^\tau\cap M^{T_2} =
M^\tau\cap M^{T},
\end{equation}
which implies (i).

Each $x\in M^T$ has a $T^\times$-equivariant neighborhood $U_x$ on which the $T^\times$-action
is conjugate to a linear action.
The three conditions are clearly equivalent for a linear action, 
so Condition $(i)$ or $(ii)$ implies $(iii)$.

We now show by contradiction that $(iii)$ implies $(ii)$.
Suppose that $(ii)$ does not hold:
that is, there exists $x\in M^{T_2}$ with $x\notin M^T$.
Let $\Phi^\xi_t$ be the gradient flow
of $\Phi^\xi$. Then $\Phi^\xi_t$ is a $T^+$-equivariant
diffeomorphism of $M$.
Thus, $\Phi^\xi_t(x)$ has the same property of $x$ but,
if $t$ is large enough, $\Phi^\xi_t(x)$
will belong to $U_x$ for some $x\in M^T$.
This contradicts $(iii)$. 
\end{proof}

\begin{Lemma}\label{Hspeqforpri}
Let $M$ be a compact symplectic manifold
equipped with a Hamiltonian action of a torus $T$.
Let $\tau$ be a smooth
anti-symplectic involution on $M$ compatible with the action of $T$.
Suppose that $M^T=M^{T_2}$ and that
$\pi_0(M^T\cap M^\tau)\to \pi_0(M^T)$ is a bijection.
Then $M^\tau$ is $T_2$-\ef\ over $\bbz_2$.
\end{Lemma}

\begin{proof}
As $\pi_0(M^T\cap M^\tau)\to \pi_0(M^T)$ is a bijection,
by \cite[Lemma\,2.1 and Theorem\,3.1]{Du},
we know that $B(M^\tau)=B(M^\tau\cap M^T)$.
By Lemma \ref{nuchi}, $M^\tau\cap M^T = (M^\tau)^{T_2}$
so $B(M^\tau)=B((M^\tau)^{T_2})$.
This implies that $M^\tau$ is $T_2$-equivariantly
formal over $\bbz_2$ (see, e.g.\,\cite[Proposition\,1.3.14]{AP}). 
\end{proof}

%%%%%%%%%%%%%%%%%%%%%%%%%%%%%%%%%%%%%%%%%%%%%%%%%%%%%%%%%%%%%%%%%%%
\subsection{The main theorems}

\begin{Theorem}\label{HamilCS}
Let $M$ be a compact symplectic manifold equipped with a
Hamiltonian action of a torus $T$ and with
a compatible smooth anti-symplectic involution
$\tau$. If $M^T$ is a \csp, then $M$ is a \csp.
\end{Theorem}

\begin{proof}
Choose a generic $\xi\in\algt$ so that $\Phi^\xi \co  M\to \bbr$
is a Morse-Bott function with critical set $M^T$.
Let $c_0<c_1 <\cdots
< c_N$ be the critical values of $\Phi^\xi$, and let $F_i =
(\Phi^\xi)^{-1}(c_i)\cap M^T$ be the critical sets. Let
$\varepsilon >0$ be less than any of the differences
$c_i-c_{i-1}$, and define $M_i = (\Phi^\xi)^{-1}((-\infty,c_i
+ \varepsilon])$.
We will prove by induction that $M_i$ is a \csp.
This is true for $i=0$ since $M_0$ is $C$-homotopy
equivalent to $F_0$, which is a \csp\ by
hypothesis. By induction, suppose that $M_{i-1}$
is a \csp.

We saw in Subsection~\ref{Hampreli} that the negative normal
bundle $\nu_i$ to $F_i$  is a $\tau$-bundle.
The pair $(M_i,M_{i-1})$ is $C$-homotopy equivalent to
the pair $(\bbd(\nu_i),\bbs (\nu_i))$. Since $F_i$ is a
\csp\ by hypothesis, the pair $(M_i,M_{i-1})$
is \cpa\ by \proref{thomcbdle}. Therefore,
$M_i$ is a \csp\ by \proref{EXPlac}.
We have thus proven that each $M_i$ is a \csp, including $M_N=M$.

\end{proof}

\begin{Remark} \rm
The proof of Theorem~\ref{HamilCS} shows that
the compactness assumption on $M$ can be replaced by the
assumptions that $M^T$ consists of finitely many connected
components, and that some generic component of the moment map
$\Phi\co  M\to \algt^*$ is proper and bounded below. That $M^T$ has
finitely many connected components ensures that $H_T^*(M)$ is a
finite rank module over $H_T^*(pt)$. That some component of the
moment map is proper and bounded below ensures that that component
of the moment map is a Morse-Bott function on $M$. Examples of this
more general situation include hypertoric manifolds (see~\cite{HH}).
\end{Remark}

Using Theorem~\ref{strcomp} and Corollary~\ref{strcomp},
we get the following corollary of Theorem~\ref{HamilCS}.

\begin{Corollary}\label{HamilCSeq}
Let $M$ be a compact symplectic manifold equipped with a
Hamiltonian action of a torus $T$ and a compatible smooth
anti-symplectic involution $\tau$.
If $M^T$ is a \csp, then $M_T$ is a \csp.
In particular, there is a ring isomorphism
$$
\bar\kappa\co  H^{2*}_T(M)\hfl{\approx}{} H^{*}_{T_2}(M^\tau). \eqno{\qed}
$$
\end{Corollary}

Finally, the same proof as for \thref{HamilCS}, 
using \proref{ThomSCC} instead of \proref{thomcbdle}, gives the following:

\begin{Theorem}\label{HamilSCC}
Let $M$ be a compact symplectic manifold equipped with a
Hamiltonian action of a torus $T$ and with
a compatible smooth anti-symplectic involution
$\tau$. If $M^T$ is a \scc, then $M$ is a \scc. \qed
\end{Theorem}

\begin{Examples}\rm
The  theorems of this subsection apply to toric manifolds
($M^T$ is discrete). They also apply to spatial polygon spaces
${\rm Pol}(a)$ of $m$ edges, with lengths $a=(a_1,\dots,a_m)$
(see, e.g.~\cite{HK}), the involution being given by a mirror reflection \cite[\S,9]{HK}. One proceeds by induction $m$ (for $m\leq 3$, ${\rm Pol}(a)$
is either empty or a point). The induction step uses that  
${\rm Pol}(a)$ generically admits compatible Hamiltonian circle action, 
called bending flows, introduced by Klyachko (\cite{Kl}, see, e.g.~\cite{HT}), for which the connected component of the fixed point set are polygon 
spaces with fewer edges \cite[Lemma 2.3]{HT}. 

Therefore, toric manifolds and polygon spaces are \scc es.
The isomorphism $\kappa$ were discovered in \cite{DJ} and \cite[\S\,9]{HK}.
\end{Examples}

%%%%%%%%%%%%%%%%%%%%%%%%%%%%%%%%%%%%%%%%%%%%%%%%%
\subsection{The Chevalley involution on co-adjoint orbits of
semi-simple compact Lie groups}\label{Cheva}\label{coorb}

The goal of this section is to show that coadjoint orbits of
compact semi-simple Lie groups are equipped with a
natural involution which makes them \csp s. Let $\algl$
be a semi-simple complex Lie algebra, and $\algh$ a Cartan
sub-algebra with roots $\Delta$. Multiplication by $-1$ on
$\Delta$ induces, by the isomorphism theorem
\cite[Corollary C,\S\,2.9]{Sa}, a Lie algebra involution $\sigma$ on
$\algl$ called the {\it Chevalley involution} \cite[Example
p.\,51]{Sa}. Then $\sigma(h)=-h$ for $h\in\algh$ and
$\sigma(X_\alpha)=-X_{-\alpha}$, if $X_\alpha$ is the weight
vector occurring in a Chevalley normal form \cite[Theorem A,\S\,2.9
]{Sa}. By construction of the compact form $\algl_0$ of $\algl$
\cite[\S\,2.10]{Sa}, the involution $\sigma$ induces a Lie algebra
involution on the real Lie algebra $\algl_0$, still called the
Chevalley involution and denoted by $\sigma$. This shows that any
semi-simple compact real Lie algebra admits a Chevalley
involution. For instance, if $\algl=\mathfrak{sl}(n,\bbc)$, then
$\sigma(X)=-X^T$ and the induced Chevalley involution on
$\algl_0=\mathfrak{su}(n)$ is complex conjugation.

Let $G$ be a compact semi-simple Lie group with Lie algebra
$\algg$ and a maximal torus $T$. Recall that the dual $\algg^*$ of
$\algg$ is endowed with a Poisson structure characterized by the
fact that $\algg^{**}$ is a Lie sub-algebra of
$\calc^\infty(\algg)$ and the canonical map
$\algg\hfl{\approx}{}\algg^{**}$ is a Lie algebra isomorphism.
Therefore, the map $\tau = -\sigma^*\co \algg^*\to\algg^*$ is an
anti-Poisson involution, called again the {\it Chevalley involution} on
$\algg^*$.

\begin{Theorem}\label{invtaucoad}
The Chevalley involution $\tau$ preserves each coadjoint orbit $\calo$, and
induces an anti-symplectic involution $\tau\co  \calo\to\calo$ with
respect to which $\calo$ is a \csp. One also has an ring-isomorphism
$$
\bar\kappa\co  H^{2*}_T(\calo)\hfl{\approx}{} H^{*}_{T_2}(\calo^\tau).
$$
\end{Theorem}

\begin{proof}
The coadjoint orbits are the symplectic leaves of the
Poisson structure on $\algg^*$. As $\tau$ is anti-Poisson, the
image $\tau(\calo)$ of a coadjoint orbit $\calo$ is also a
coadjoint orbit $\calo'$. We will show that $\calo'=\calo$. Since
$G$ is semi-simple, the Killing form $\llangle{}{}$ is negative
definite. Thus, the map $K\co \algg\to\algg^*$ given by
$K(x)(-)=\llangle{x}{-}$ is an isomorphism. It intertwines the
adjoint action with the coadjoint action and satisfies $\tau\pcirc
K=-K\pcirc\sigma$.

Now, to show $\calo'=\calo$, if $\calo$ is a coadjoint orbit, the
adjoint orbit $K^\mun(\calo)$ contains an element $t\in\algt$.
Thus, $\tau(K(t))= -K(\sigma(t))=K(t)$. Therefore
$\calo'=\tau(\calo)=\calo$. As $\tau$ is anti-Poisson on
$\algg^*$, its restriction to $\calo$ is anti-symplectic.
Moreover, since $\sigma$ is $-1$ on $\algt$, the involution $\tau$
is compatible with the coadjoint action of $T$ on $\calo$.
Finally, $\calo^T$ is discrete, and  $\calo\cap
K(\algt)=\calo^T\subset\calo^\tau$.  It is clear, then, that
$\calo^T$ is a \csp.  The theorem now follows from
Theorem~\ref{HamilCS} and Corollary~\ref{HamilCSeq}.
\end{proof}

\begin{Remark}
The conjugation cells used to build $\calo$ as a conjugation space
are precisely the Bruhat cells of the coadjoint orbit.  The Bruhat
decomposition is $\tau$-invariant.
\end{Remark}

\begin{figure}[ht!]\small\anchor{fig}
%\centerline{
\begin{center}
\psfrag{1}{(a)}
\psfrag{2}{\! (b)}\psfrag{3}{\! (c)}\psfrag{a}{$\alpha$} \psfrag{b}{$\beta$}
\psfrag{b2a}{$\beta + 2\alpha$}
\epsfig{figure=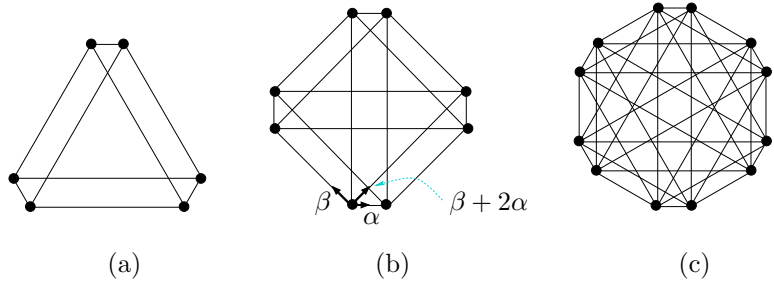,width=4in}
%} \vskip 0mm \centerline{
\parbox{4.8in}{\caption{The moment polytopes
for the generic coadjoint orbits of simple Lie groups of rank $2$:
we show types  (a) $A_2$, (b) $B_2$ and (c) $G_2$. As shown in (b), for type $B_2$,
at a $T$-fixed point, we can see that $\beta$, $\alpha$ and $\beta+
2\alpha$ are isotropy weights.  There is a similar occurrence for
type $G_2$. }\label{fig}}
%}
\end{center}
\end{figure}

In type $A$, the Chevalley involution is complex conjugation on
$\mathfrak{su}(n)$.  In this case, Theorem~\ref{invtaucoad} has
been proven in \cite{Sc} and \cite{BGH}. In those papers, the
authors use the fact that the isotropy weights at each fixed point
are pairwise independent over $\bbf_2$. This condition is not
satisfied in general for the coadjoint orbits of other types.
Indeed, for the generic orbits, these weights are a set of
positive roots and the other types have strings of roots of length
at least $2$. This can be seen already in the moment polytopes for
generic coadjoint orbits of $B_2$ and $G_2$, shown in \figref{fig}.

In \cite{Sc} and \cite{BGH}, the isomorphism
$$
\bar\kappa\co  H^{2*}_T(\calo)\hfl{\approx}{} H^{*}_{T_2}(\calo^\tau)
$$
is proved by giving a combinatorial description of each of these
rings, and noting that these descriptions are identical.  This
combinatorial description does not generally apply in the other
types precisely because the isotropy weights the fixed points are
not pairwise independent over $\bbf_2$.  Nevertheless, we still
have the isomorphism on the equivariant cohomology rings.

%%%%%%%%%%%%%%%%%%%%%%%%%%%%%%%%%%%%%%%%%%%%%%%%%
\subsection{Symplectic reductions}\label{sred}

Let $M$ be a compact symplectic manifold equipped with a Hamiltonian
action of a torus $T$ and a compatible smooth anti-symplectic
involution $\tau$. We saw in Theorem~\ref{HamilCS} that if $M^T$ is a
\csp, then $M$ is a \csp.  Using this, we extend results of Goldin and
the second author \cite{GH} to show that in certain cases, the
symplectic reduction is again a \csp.  To do this, we must construct a
ring isomorphism
$$
\ri_{red}\co H^{2*}(M/\!/T(\mu))\to H^{*}((M/\!/T(\mu))^{\tau_{red}}),
$$
and a section
$$
\sigma_{red}  \co  H^{2*}(M/\!/T(\mu)) \to H^{2*}_C(M/\!/T(\mu))
$$
that satisfy the conjugation equation.

Let $ \Phi  \co  M \to  \algt^* $  be  the moment map for $ M $. When
$\mu\in\algt^*$ is a  regular value of $\Phi$, and when $T$ acts
on $\Phi^{-1}(\mu)$ freely, we define the symplectic reduction
$$
M/\! /T(\mu) = \Phi^{-1}(\mu)/T.
$$
Kirwan \cite{kirwan} proved that the inclusion map
$\Phi^{-1}(\mu)\into M$ induces a surjection in equivariant cohomology
with rational coefficients:
\begin{equation}\label{eq:surjection}
\xymatrix{
H_T^*(M;\bbq)\ar@{->>}[r]^(0.25){\calk} &
H_T^*(\Phi^{-1}(\mu);\bbq)=H^*(M/\! /T(\mu);\bbq).
}
\end{equation}
The map $\calk$ is called the Kirwan map.  Under additional
assumptions on the torsion of the fixed point sets  and the group
action, this map is surjective over the integers or $\bbz_2$ as
well.  There are several ways to compute the kernel of $\calk$.
Tolman and Weitsman \cite{TW:ker} did so in the way that
is most suited to our needs.

Goldin and the second author  extend these two results to the real
locus, when the the torus action has suitable $2$-torsion.
\begin{Definition}\label{twotordefi}
Let $x\in M$, and suppose $H$ is the identity component of the
stabilizer of $x$. Then we say $x$ is a {\it $2$-torsion point} if
there is a weight $\alpha$ of the isotropy action of $H$ on the normal
bundle $\nu_xM^H$ that satisfies $\alpha\equiv 0\ \mod\ 2$.
\end{Definition}
The necessary assumption is that $M^\tau$ have no $2$-torsion points.
This hypothesis is reasonably  strong. Real loci of  toric varieties and
coadjoint orbits in type  $A_n$ satisfy this hypothesis, for
example, but the real loci of  maximal  coadjoint orbits  in type
$B_2$  do not.

We now  define reduction in the  context of real  loci.  Fix
$\mu$ a  regular  value  of $\Phi$ satisfying the condition that $T$
acts freely on $\Phi^{-1}(\mu)$.  Then $M_{red}=M/\!/T(\mu)$ is again
a symplectic manifold with a canonical  symplectic form
$\omega_{red}$.  Moreover,  there is an induced involution
$\tau_{red}$  on $M_{red}$, and this  involution is
anti-symplectic.  Thus, the fixed point set of this involution
$
(M/\!/T(\mu))^{\tau_{red}}
$
is a Lagrangian submanifold of $M$.  We now define
$$
M^\tau/\!/T_2(\mu)\ccoeg  ((\Phi|_{M^\tau})^{-1}(\mu))/T_2.
$$
When $T$ acts freely on the level set, Goldin and the second
author \cite{GH} show that
$$
(M/\!/T(\mu))^{\tau_{red}}=M^\tau/\!/T_2(\mu).
$$
We can now start proving that, under certain
hypotheses, the quotient $M/\!/T(\mu)$ is a conjugation space. We
begin by constructing the isomorphism $\ri_{red}$.

\begin{Proposition}\label{prokappared}
Suppose $M$ is a compact symplectic manifold equipped with a
Hamiltonian action of a torus $T$ and a compatible smooth
anti-symplectic involution $\tau$. Suppose further that $M^T$ is a
\csp,  and that $M$ contains no $2$-torsion points. Then there is
an isomorphism
$$
\ri_{red}\co H^{2*}(M/\!/T(\mu))\hfl{\approx}{}H^{*}(M^\tau/\!/T(\mu)) =
H^{*}((M/\!/T(\mu))^{\tau_{red}}),
$$
induced by $\ri$.
\end{Proposition}

\begin{proof}
The first main theorem of \cite{GH}  states that when $M^\tau$
contains no $2$-torsion  points,  the   {\em real   Kirwan   map}
in equivariant cohomology
$$
\calk^\tau\co H^*_{T_2}(M^\tau)\to
H^*_{T_2}(\Phi|_{M^\tau}^{-1}(\mu))=H^*(M^\tau/\! /T(\mu)),
$$
induced by inclusion, is a surjection.  The proof of surjectivity
makes use of the function $|\!|\Phi-\mu|\!|^2$ as a Morse-Kirwan
function on $M^\tau$.  The critical sets of this function are possibly
singular, but the hypothesis that the real locus have no $2$-torsion
points allows enough control over these critical sets to prove
surjectivity.

Let $x\in M^T$. By assumption $x$ is not a $2$-torsion point, so
Condition~\eqref{nuchi3} of Lemma~\ref{nuchi} is
satisfied. Lemma~\ref{nuchi} then implies that $M^T=M^{T_2}$. We now
show that there is a commutative diagram
\begin{equation}\label{eq:diagramInj}
\begin{tabular}{c}
\xymatrix@C+3pt@M+4pt@R-4pt{%
 H_T^{2*}(M) \ar[d]_{\bar\ri}^{\iso} \ar@{>->}[r]^{i} &
H_T^{2*}(M^T) \ar[d]^{\iso}_{\ri} \\
H_{T_2}^{*}(M^\tau) \ar@{>->}[r]^{i^\tau} &
H_{T_2}^{*}((M^\tau)^{T_2}) }
\end{tabular}
\end{equation}
where the horizontal arrows are induced by inclusions.  To see this,
we first note that because $M^T$ is a \csp, then $M_T$ is a \csp\
by Corollary~\ref{HamilCSeq}, which also gives
the left isomorphism $\bar\ri$. The
trivial $T$-action on $M^T$ is also compatible with $\tau$.
By Theorem~\ref{strcomp}, one have a ring isomorphism
$\kappa\co H_T^{2*}(M^T)\hfl{\approx}{}H_{T_2}^*(M^\tau\cap M^T)$. As
$M^T=M^{T_2}$, we deduce that $M^\tau\cap M^T=(M^\tau)^{T_2}$
by Lemma~\ref{nuchi}, whence the the right vertical isomorphism $\ri$.

Diagram~\eqref{eq:diagramInj} is commutative by the naturality of
\hfra s (Proposition \ref{eqmapcoh}).
Finally, $M^\tau$ is $T_2$-\ef\ over $\bbz_2$  by Lemma~\ref{Hspeqforpri}.
Therefore $i^\tau$ is injective by,
e.g.\,\cite[Proposition\,1.3.14]{AP}. It follows that $i$ is also
injective.

Note that Kirwan showed that $i$ is injective when the coefficient ring is
$\bbq$.  However, an additional assumption on $M^T$ is needed to extend her
proof to the coefficient ring $\bbz_2$, so we may not conclude that
directly.

We denote the restriction  of a  class  $\alpha\in H_{T}^{*}(M^\tau)$
to the fixed points by $\alpha|_{(M^\tau)^{T_2}}\in
H_{T_2}^{*}((M^\tau)^{T_2})$.
The  second   main  result  of   \cite{GH}  computes  the   kernel
of $\calk^\tau$.  For every $\xi\in\algt$, let
$$
M^\tau_\xi  = \left\{ p\in M^\tau\  |\ \left<
\Phi(p),\xi\right>\leq 0\right\} \subseteq M^\tau.
$$
Let $F=M^T$ denote the fixed point set, and let
$$
K^\tau_\xi = \bigg\{ \alpha\in H_{T}^*(M^\tau)\ \bigg|\
\alpha|_{F\cap M^\tau_\xi} = 0\bigg\}.
$$
Finally, let $K^\tau$ be the ideal generated by the ideals
$K^\tau_\xi$ for all $\xi\in \algt$.  Then there   is  a short
exact sequence,  in cohomology  with $\bbz_2$ coefficients,
$$
0\to K^\tau\to H^*_{T}(M^\tau)\to H^*(M^\tau/\! /T(\mu))\to 0.
$$
The important thing  to notice is that this  description of the
kernel is {\em identical} to the description  of the kernel for
$M$, given by Tolman  and Weitsman, when  $M$ contains  no
$2$-torsion  points.  The fact  that  Diagram~\ref{eq:diagramInj}
commutes  implies  that  the support of a  class $\ri(\alpha)$ is
the real locus  of the support of $\alpha$.  Therefore,  there is
a natural isomorphism  between $K$ and $K^\tau$ induced by $\ri$.
Thus, we have a commutative diagram:
$$
\xymatrix{ 0\ar[r] & K \ar[r]\ar[d]_{\iso} &
H_T^{2*}(M)\ar[r]^<<<<<{\calk}\ar[d]_{\ri}^{\iso} &
    H^{2*}(M/\!/T(\mu))\ar[r]\ar@{-->}[d] & 0 \\
0\ar[r] & K^\tau \ar[r] &
H_{T}^{*}(M^\tau)\ar[r]^<<<<<{\calk^\tau} &
    H^{*}(M^\tau/\!/T(\mu))\ar[r] & 0}
$$
Therefore, the vertical dashed arrow represents an induced isomorphism
\begin{equation}
\ri_{red}\co H^{2*}(M/\!/T(\mu))\hfl{\approx}{} H^{*}(M^\tau/\!/
T(\mu)),
\end{equation}
as rings.
\end{proof}

Now that we have established the isomorphism $\ri_{red}$ between
the cohomology of the symplectic reduction and the cohomology of
its real points, we must find the map $\sigma_{red}$ and prove the
conjugation relation.  We have the following commutative diagram:
$$\xymatrix{
H^{2*}(M_T)  \ar@{>>}[d]_{\mathcal{K}}
\ar@/^/[r]^{\sigma} &
H_C^{2*}(M_T) \ar[d]_{\mathcal{K}_C} \ar@{>>}[l]^\rho \\
H^{2*}(M/\! /T)
 & H_C^{2*}(M/\! /T)
\ar[l]_{\rho_{\mathrm{red}}}}
$$
As the diagram commutes, we see that $\rho_{\mathrm{red}}$ is a surjection.
Moreover, because $\mathcal{K}$ is a surjection, we may choose an
additive section $s \co  H^{2*}(M/\! /T) \to H^{2*}(M_T)$ and then
define a section $\sigma_{\mathrm{red}} \ccoeg  \mathcal{K}_C\circ \sigma\circ
s$ of $\rho_{\mathrm{red}}$. Adding the restriction maps into the diagram,
we have:
\begin{equation}\label{eq:commdiagr}
 \xymatrix{ H^{2*}(M_T)
\ar@{>>}[d]_{\mathcal{K}} \ar@/^/[r]^{\sigma} &
H_C^{2*}(M_T)
\ar[d]_{\mathcal{K}_C} \ar@{>>}[l]^\rho\ar[r]^<<<<<<<r &
H_C^{2*}(M_T^\tau) \iso H^{2*}(M_T^\tau)[u]
 \ar@<5ex>[d]_{\mathcal{K}^\tau\otimes 1} \\
 H^{2*}(M/\! /T) \ar@/_/[u]_{s}
\ar@/_/[r]_{\sigma_{\mathrm{red}}} & H_C^{2*}(M/\! /T)
\ar@{>>}[l]_{\rho_{\mathrm{red}}} \ar[r]_(0.3){r_{\mathrm{red}}} &
H_C^{2*}((M/\! /T)^\tau) \iso H^{*}((M/\! /T)^\tau)[u]
 }
\end{equation}
Now we check, for $a\in H^{2m}(M/\! /T)$,
\begin{eqnarray*}
r_{red}(\sigma_{red}(a)) & = & r_{\mathrm{red}} (\mathcal{K}_C\circ
\sigma\circ s (a))\\
& = & \mathcal{K}^\tau\otimes 1 ( r(\sigma ( s(a)))) \\
& = & \mathcal{K}^\tau\otimes 1 ( \kappa(s(a)) u^m + \lt{m})\\
& = & \kappa_{\mathrm{red}}(a)u^m + \lt{m}.
\end{eqnarray*}
Thus, by the commutativity of diagram \eqref{eq:commdiagr}, we have
proved the conjugation equation, and hence the following theorem.

\begin{Theorem}\label{Thsympred}
Let $M$ be compact symplectic manifold equipped with a
Hamiltonian action of a torus $T$, with moment map $\Phi$,
and with a compatible smooth anti-symplectic involution $\tau$.
Suppose that $M^T$ is a \csp\ and that $M$ contains no $2$-torsion points.
Let $\mu$ be a regular value of $\Phi$ such that $T$ acts
freely on $\Phi^\mun(\mu)$. Then, $M/\! /T(\mu)$
is a conjugation space. \qed
\end{Theorem}

\begin{Remark}\label{remcuts}\rm
When $T=S^1$ in Theorem~\ref{Thsympred}, the symplectic cuts $C_\pm$ at $\mu$
introduced by E.~Lerman \cite{Le} also inherit an Hamiltonian
$S^1$-action and a compatible anti-symplectic involution.
The connected components of $C_\pm^T$ are those of $M^T$ plus
a copy of $M/\! /T(\mu)$. By Theorem~\ref{Thsympred}, $C_\pm^T$
are \csp s. Therefore, using Theorems~\ref{HamilCS} and Corollary~\ref{HamilCSeq},
we deduce that $C_\pm$ and $(C_\pm)_T$ are \csp s.
\end{Remark}

\Addresses\recd

\end{document}

%% file: agtout.tex
\def\ifplaintex{\expandafter\ifx\csname documentclass\endcsname\relax}

\def\gtp{{\mathsurround=0pt\it $\cal G\mskip-2mu$eometry \&\ 
$\cal T\!\!$opology $\cal P\!$ublications}}  % GT publications

\def\recd{{\small Received:\qua\receiveddate\ifx\reviseddate\relax
\else\qquad Revised:\qua\reviseddate\fi\par}} 

\def\lognumber#1{\def\thelognumber{#1}}
\def\volumenumber#1{\def\thevolumenumber{#1}}
\def\volumeyear#1{\def\thevolumeyear{#1}}
\def\papernumber#1{\def\thepapernumber{#1}}
\def\pagenumbers#1#2{\def\startpage{#1}\def\finishpage{#2}}
\def\published#1{\def\publishdate{#1}}

\def\received#1{\def\receiveddate{#1}}

\def\accepted#1{\def\accepteddate{#1}}

\def\asciiauthors#1{\def\theasciiauthors{#1}}
\def\asciiaddress#1{\def\theasciiaddress{#1}}
\def\asciiemail#1{\def\theasciiemail{#1}}

\long\def\asciiabstract#1{\long\def\theasciiabstract{#1}}

\let\\\par\let\thelognumber\relax\let\thevolumenumber\relax
\let\thepapernumber\relax\let\thevolumeyear\relax\let\startpage\relax
\let\finishpage\relax\let\publishdate\relax\let\receiveddate\relax
\let\reviseddate\relax\let\accepteddate\relax\let\theasciititle\relax
\let\theasciiauthors\relax\let\theasciiaddress\relax
\let\theasciiabstract\relax

\let\theasciiemail\relax

\ifplaintex
\font\logobig=cmssbx10 scaled 3836
\font\logomed=cmssbx10 scaled 2557
\else
\font\logobig=cmssbx10 scaled 4200
\font\logomed=cmssbx10 scaled 2800
\fi

\long\def\makeagttitle{   %%% start of definition of \makeagttitle
\count0=\startpage
\agt\hfill      %   Journal title (top left) 
\hbox to 45truept{\vbox to 0pt{\vglue -13truept{\logomed A\kern -.37em{\logobig 
T}\kern -.38em G}\vss}\hss}
\break
{\small Volume \thevolumenumber\ (\thevolumeyear)
\startpage--\finishpage\nl
Published: \publishdate}

\vglue .25truein

{\parskip=0pt\leftskip 0pt plus
1fil\def\\{\par\smallskip}{\Large\bf\thetitle}\par\medskip} \vglue
0.05truein

{\parskip=0pt\leftskip 0pt plus 1fil\def\\{\par}{\sc\theauthors}
\par\medskip}%
 
\vglue 0.03truein 

{\small\leftskip 25truept\rightskip 25truept{\bf Abstract}\stdspace\theabstract

{\bf AMS Classification}\stdspace\theprimaryclass
\ifx\thesecondaryclass\relax\else; \thesecondaryclass\fi\par
{\bf Keywords}\stdspace \thekeywords\par}\vglue 7truept

}   %%%% end of definition of \makeagttitle

\ifplaintex
\font\phead=cmsl9 scaled 950
\font\pnum=cmbx10 scaled 913
\font\pfoot=cmsl9 scaled 950
\headline{\vbox to 0pt{\vskip -4.5mm\line{\small\phead\ifnum
\count0=\startpage ISSN 1472-2739 (on-line) 1472-2747 (printed)
\hfill {\pnum\folio}\else\ifodd\count0\def\\{ }% 
\ifx\theshorttitle\relax\thetitle\else\theshorttitle\fi\hfill{\pnum\folio}
\else\def\\{ and }{\pnum\folio}\hfill\ifx\theshortauthors\relax\theauthors
\else\theshortauthors\fi\fi\fi}\vss}}
\footline{\vbox to 0pt{\vglue 0mm\line{\small\pfoot\ifnum\count0=\startpage
\copyright\ \gtp\hfill\else
\agt, Volume \thevolumenumber\ (\thevolumeyear)\hfill\fi}\vss}}
\else
\font=cmsl9 scaled 1050
\font\lnum=cmbx10 
\font=cmsl9 scaled 1050
\makeatletter
\def\@oddhead{{\small\lhead\ifnum\count0=\startpage ISSN 1472-2739 
(on-line) 1472-2747 (printed)\hfill {\lnum\number\count0}\else\ifodd\count0
\def\\{ }\ifx\theshorttitle\relax \thetitle \else\theshorttitle\fi\hfill
{\lnum\number\count0}\else\def\\{ and }{\lnum\number\count0}
\hfill\ifx\theshortauthors\relax 
\theauthors\else\theshortauthors\fi\fi\fi}}\def\@evenhead{\@oddhead}
\def\@oddfoot{\small\lfoot\ifnum\count0=\startpage\copyright\ \gtp\hfill\else
\agt, Volume \thevolumenumber\ (\thevolumeyear)\hfill\fi}
\def\@evenfoot{\@oddfoot}
\makeatother
\fi
\let\maketitlepage\makeagttitle
\let\makeshorttitle\maketitlepage
\let\maketitle\maketitlepage

\newwrite\gtoutfile
\long\gdef\makeheadfile{  %%% start of definition of \makeheadfile
{\def\\{, }\def\s{ }
\immediate\openout\gtoutfile head.xxx
\immediate\write\gtoutfile{Proxy-for: \ifx\theasciiauthors\relax
\theauthors\else\theasciiauthors\fi\s<\ifx\theasciiemail\relax\theemail\else\theasciiemail\fi>}
\immediate\write\gtoutfile{\noexpand\\}
\immediate\write\gtoutfile{Authors: \ifx\theasciiauthors\relax
\theauthors\else\theasciiauthors\fi}
{\def\\{ }\immediate\write\gtoutfile{Title: \ifx\theasciititle\relax
\thetitle\else\theasciititle\fi}}
\immediate\write\gtoutfile{Subj-class: GT or SG, GR etc}
\immediate\write\gtoutfile{MSC-class: \theprimaryclass\ifx\thesecondaryclass\relax\else, \thesecondaryclass\fi}
\immediate\write\gtoutfile{Journal-ref: Algebr. Geom. Topol. \thevolumenumber\s
(\thevolumeyear) \startpage-\finishpage}
\immediate\write\gtoutfile{Comments: Published by Algebraic and
Geometric Topology at}
\immediate\write\gtoutfile{\s\s\s  http://www.maths.warwick.ac.uk/agt/AGTVol\thevolumenumber/agt-\thevolumenumber-\thepapernumber.abs.html}
\immediate\write\gtoutfile{\noexpand\\}
\immediate\write\gtoutfile{}
\ifx\theasciiabstract\relax
\immediate\write\gtoutfile{\theabstract}\else
\immediate\write\gtoutfile{\theasciiabstract}\fi
\immediate\write\gtoutfile{}
\immediate\write\gtoutfile{\noexpand\\}
\immediate\write\gtoutfile{}
\immediate\closeout\gtoutfile}}  %%% end of definition of \makeheadfile

\def\maketitlepage{\makeagttitle\makeheadfile}
\let\makeshorttitle\maketitlepage
\let\maketitle\maketitlepage